\RequirePackage{fix-cm}
\documentclass[smallextended,natbib]{svjour3}       
\smartqed  

\usepackage{graphicx}
\usepackage[top=2cm,bottom=2cm,left=2cm,right=2cm]{geometry}

\usepackage{amssymb}
\usepackage{siunitx}
\usepackage[colorlinks = true,
linkcolor = blue,
urlcolor  = blue,
citecolor = blue,
anchorcolor = blue]{hyperref}
\usepackage{physics} 
\usepackage{bm} 

\makeatletter
\def\cl@chapter{\@elt {theorem}}
\makeatother
\usepackage[capitalise]{cleveref}

\crefname{equation}{}{}

\usepackage{tikz}
\usepackage{svg}
\usetikzlibrary{arrows,arrows.meta,calc,decorations.markings,decorations.markings,math,patterns,positioning,shadows.blur,shapes,shapes.symbols,spy}
\usepackage{subcaption}
\usepackage{algorithmic}
\usepackage{xurl}
\usepackage{svg}
\usepackage{adjustbox}
\usepackage{pgfplots}
\usepackage{tikz}
\usepackage{framed} 
\usetikzlibrary{arrows.meta}
\usepackage{subcaption}
\usepackage{pdfpages}

\newcommand{\mm}[1]{\textcolor{teal}{#1}}

\journalname{Optimization and Engineering}
\begin{document}

\title{Combined Parameter and Shape Optimization of Electric Machines with Isogeometric Analysis
}


\author{Michael Wiesheu         \and
        Theodor Komann  \and
        Melina Merkel \and
        Sebastian Schöps \and 
        Stefan Ulbrich \and
        Idoia {Cortes Garcia} 
}


\institute{M. Wiesheu \at
              Schloßgartenstraße~8 \\
              64289 Darmstadt \\
              Germany \\
              \email{michael.wiesheu@tu-darmstadt.de}
}

\date{Received: date / Accepted: date}

\maketitle

\begin{abstract}
In structural optimization, both parameters and shape are relevant for the model performance. Yet, conventional optimization techniques usually consider either parameters or the shape separately. This work addresses this problem by proposing a simple yet powerful approach to combine parameter and shape optimization in a framework using Isogeometric Analysis (IGA). The optimization employs sensitivity analysis by determining the gradients of an objective function with respect to parameters and control points that represent the geometry. The gradients with respect to the control points are calculated in an analytical way using the adjoint method, which enables straightforward shape optimization by altering these control points. Given that a change in a single geometry parameter corresponds to modifications in multiple control points, the chain rule is employed to obtain the gradient with respect to the parameters in an efficient semi-analytical way.
The presented method is exemplarily applied to nonlinear 2D magnetostatic simulations featuring a permanent magnet synchronous motor and compared to designs, which were optimized using parameter and shape optimization separately. It is numerically shown that the permanent magnet mass can be reduced and the torque ripple can be eliminated almost completely by simultaneously adjusting rotor parameters and shape. The approach allows for novel designs to be created with the potential to reduce the optimization time substantially.




\end{abstract}


\newcommand{\NurbsBasis}{\hat{N}}
\newcommand{\BSplineBasis}{\hat{B}}
\newcommand{\NurbsDegree}{p}
\newcommand{\NumberBasisFunctions}{n}
\newcommand{\NurbsWeight}{w}
\newcommand{\Jac}{\mathbf{J}_F}
\newcommand{\GeometryBasis}{G}

\newcommand\st{\mathrm{st}}
\newcommand\rt{\mathrm{rt}}
\newcommand\ag{\mathrm{ag}}
\newcommand\source{\mathrm{src}}

\newcommand{\der}{\operatorname{d}}
\newcommand{\MagVecPot}{\mathbf{A}}
\newcommand{\MagVecPotz}{A_z}
\newcommand{\MagVecPotzRt}{A_{z,\rt}}
\newcommand{\MagVecPotzSt}{A_{z,\st}}
\newcommand{\CurrentDensity}{\mathbf{J}}
\newcommand{\CurrentDensityTwoD}{J_{z,\source}}

\newcommand{\MagFluxDensity}{\mathbf{B}}
\newcommand{\Remanence}{\mathbf{B_\mathrm{r}}}
\newcommand{\RemanenceRed}{\mathbf{B^\bot_\mathrm{r}}}
\newcommand{\RotAngle}{\beta}
\newcommand{\MagnetAngle}{\alpha}
\newcommand{\CurrentAngle}{\varphi_0}


\newcommand{\ansatzFunction}{N} 
\newcommand{\ansatzFunctionIndex}{j}
\newcommand{\testSymbol}{v}
\newcommand{\testFunction}{N} 
\newcommand{\testFunctionIndex}{i}

\newcommand{\couplingSymbol}{G}
\newcommand{\couplingMatrix}{\mathbf{\couplingSymbol}}
\newcommand{\stiffnessSymbol}{K}
\newcommand{\stiffnessMatrix}{\mathbf{\stiffnessSymbol}}

\newcommand{\solutionSymbol}{u}
\newcommand{\solutionVector}{\mathbf{\solutionSymbol}}
\newcommand{\rhsSymbol}{b}
\newcommand{\rhsVector}{\mathbf{\rhsSymbol}}
\newcommand{\optiSymbol}{x}
\newcommand{\optiVector}{\mathbf{\optiSymbol}}
\newcommand{\stateSymbol}{e}
\newcommand{\stateVector}{\mathbf{\stateSymbol}}

\newcommand{\mortarSymbol}{\lambda}
\newcommand{\mortarVector}{\bm{\mortarSymbol}}

\newcommand{\opt}{\mathrm{opt}}
\newcommand{\fopt}{f_\opt}
\newcommand{\RotMat}{\mathbf{R}_\RotAngle}
\newcommand{\ControlPointSymbol}{C}
\newcommand{\OptiCtrlPoint}{\mathbf{\ControlPointSymbol}}
\newcommand{\ParameterSymbol}{P}
\newcommand{\OptiParameter}{\mathbf{\ParameterSymbol}}
\newcommand{\adjointSymbol}{\gamma}
\newcommand{\adjointVector}{\bm{\adjointSymbol}}
\newcommand{\Torque}{T}
\newcommand{\MeanTorque}{\overline{\Torque}}
\newcommand{\StdTorque}{\hat{\Torque}}
\newcommand{\Amagnet}{A_\mathrm{Magnet}}
\newcommand{\Ttarget}{\Torque_\mathrm{Target}}

\newcommand{\OmegaRotor}{\Omega_\rt}
\newcommand{\OmegaStator}{\Omega_\st}
\newcommand{\GammaAirGap}{\Gamma_\ag}
\newcommand\IntegG{\mathrm{d}\Gamma}
\newcommand\Integ{\mathrm{d}\Omega}

\newcommand{\Iapp}{I_{\mathrm{app}}}
\newcommand{\nWind}{n_{\mathrm{wind}}}
\newcommand{\Acoil}{A_{\mathrm{coil}}}
\newcommand{\polePair}{p}

\definecolor{TUDa-2a}{HTML}{009CDA}
\definecolor{TUDa-2b}{HTML}{0083CC}
\definecolor{TUDa-3d}{HTML}{0071F3}
\definecolor{TUDa-3a}{HTML}{50B695}
\definecolor{TUDa-9b}{HTML}{E6001A}

\begin{table*}[!ht]   
\begin{framed}
\textbf{Nomenclature} \\

\begin{minipage}{0.5\linewidth}
    \textit{Abbreviations}\\
     \begin{tabular}{l l}
         CAD & Computer Aided Design \\
         EA & Evolutionary Algorithm \\
         FE & Finite Element  \\
         IGA & Isogeometric Analysis \\
         NURBS & Non-Uniform Rational B-Spline \\
         PDE & Partial Differential Equation \\
     \end{tabular}\\
     
     \textit{Roman symbols} \\
     \begin{tabular}{l l}
         $\mathbf{A}$ & Magnetic vector potential (\unit{V.s.m^{-1}})\\
         $A$ & Area (\unit{m^{2}}) \\
         $\mathbf{B}$ & Magnetic flux density (\unit{T}) \\
         $\mathbf{b}$ & Magnetic force vector\\
         $\OptiCtrlPoint$ & Geometry control points \\
         $\mathbf{F}$ & Mapping \\
         $\mathbf{G}$ & Coupling Matrix \\
         $\mathbf{H}$ & Magnetic field strength (\unit{A.m^{-1}}) \\
         $\mathbf{J}$ & Electric current density (\unit{A.m^{-2}}) \\
         $\mathbf{K}$ & Stiffness Matrix \\
         $\mathbf{n}$ & Surface normal vector \\
         $N$ & Basis function \\
         $\OptiParameter$ & Geometry parameters \\
         $p$ & Polynomial degree \\
        
         $\RotMat$ & Rotation matrix \\
        
         $\mathbf{u}$ & Discretized magnetic vector potential\\
         $\NurbsWeight$ & Weight \\
     \end{tabular}

\end{minipage}
\begin{minipage}{0.5\linewidth}

    \textit{Greek symbols}\\
     \begin{tabular}{l l}
         $\beta$ & Rotation angle \\
         $\Gamma$ & Domain surface \\
         $\bm{\lambda}$ & Lagrangian multiplier \\
         $\bm{\mu}$ & Magnetic permeability (\unit{H.m^{-1}}) \\
         $\bm{\nu}$ & Magnetic reluctivity (\unit{m.H^{-1}}) \\
         $\Omega$ &  Physical domain \\
         $\hat{\Omega}$ & Parametric domain \\
         $\CurrentAngle$ & Electric phase offset \\
         $\theta$ & Evaluation angle in the air gap \\
         $\Xi$ & Knot vector \\
         $\xi$ & Knot vector entry \\
     \end{tabular}\\
    
    \textit{Subscripts}\\
     \begin{tabular}{l l}
         $\ag$ & air gap \\
         azi & azimuthal \\
         max & maximum \\
         $\opt$ & optimization \\
         pm & permanent magnet \\
         r & remanent \\
         rad & radial \\
         $\rt$ & rotor \\
         $\st$ & stator \\
     \end{tabular}\\
     
\end{minipage}


\end{framed}
\end{table*}

\section{Introduction}
\label{sec:Introduction}

With the continuing trend towards electrification in all major sectors, including the industrial, residential and transportation ones,
electric motors are steadily increasing in importance. 
One example illustrating this transformation can be found in the automotive industry: To achieve the net-zero scenario by 2050 \citep{IEA_2021xx}, the International Energy Agency projects that the share of electric vehicle sales needs to increase from around \SI{10}{\percent} in 2021 to almost \SI{60}{\percent} in 2030 \citep{IEA_2022xx}. The use of computational methods is critical to generate better machine designs and enable this rapid development. Having accurate virtual prototypes greatly reduces development time, especially in the process of optimization, when a large number of different geometry configurations need to be evaluated in order to generate optimal designs. A common way to improve the performance of a design is to carry out structural optimization by making changes to its geometry.

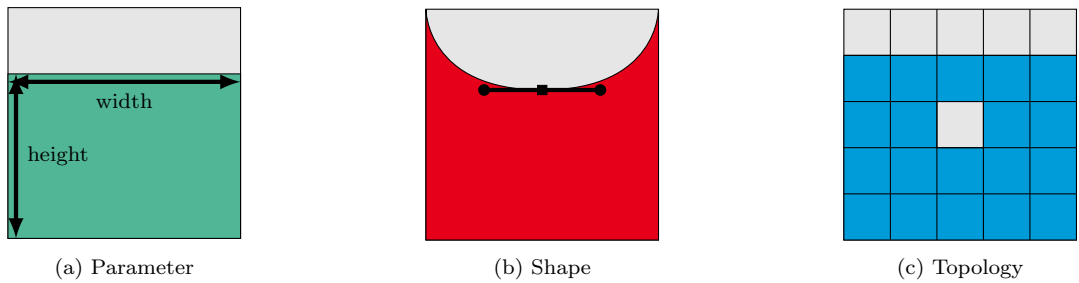
\begin{figure}[h]
    \centering
    \newcommand{\tw}{0.6\linewidth}
    \begin{subfigure}{.3\linewidth}
    \centering
    \begin{tikzpicture}[y=1pt, x=1pt, yscale=-1, xscale=1, inner sep=0pt, outer sep=0pt]
					\draw[fill=black!10] (  0.0,  0.0) rectangle (\tw,\tw);
				    \draw[fill=TUDa-3a] (  0.0, 25.0) rectangle (\tw,\tw);
					\draw[latex-latex, ultra thick]   (  3.0, 25.0) -- (  3.0,\tw) node [midway, right=0.5em] {height};
					\draw[latex-latex, ultra thick]   (  0.0, 28.0) -- (\tw, 28.0) node [midway, below=0.5em] {width};
				\end{tikzpicture}
    \caption{Parameter}
    \label{fig:Optimization:Parameter}
    \end{subfigure} \hspace{0.2cm}
    \begin{subfigure}{.3\linewidth}
    \centering
      \begin{tikzpicture}[y=1pt, x=1pt, yscale=-1, xscale=1, inner sep=0pt, outer sep=0pt]
			\draw[fill=black!10] (  0.0,  0.0) rectangle (\tw,\tw);
			\draw[fill=TUDa-9b] (  0.0,\tw) -- (  0.0,  0.0) .. controls (  0.0,0.0000) and (  0.0,0.35*\tw) .. ( 0.5*\tw,0.35*\tw) .. controls (\tw,0.35*\tw) and (\tw,0.0000) .. (\tw,  0.0) -- (\tw,\tw) -- cycle;
			\draw[fill=black] ( 0.25*\tw, 0.35*\tw) circle (2);
			\draw[fill=black] ( 0.75*\tw, 0.35*\tw) circle (2);
			\draw[fill=black] (0.48*\tw,0.33*\tw) rectangle (0.52*\tw,0.37*\tw);
			\draw[ultra thick,black] ( 0.25*\tw, 0.35*\tw) -- ( 0.75*\tw, 0.35*\tw);
		\end{tikzpicture}
		\caption{Shape}
		\label{fig:Optimization:Shape}
    \end{subfigure} \hspace{0.2cm}
    \begin{subfigure}{.3\linewidth}
    \centering
      \begin{tikzpicture}[y=1pt, x=1pt, yscale=-1, xscale=1, inner sep=0pt, outer sep=0pt]
			\draw[fill=black!10] (  0.0,  0.0) rectangle (\tw,\tw);
			\draw[fill=TUDa-2a ] (  0.0, 0.2*\tw) rectangle (\tw,\tw);
			\draw[fill=black!10] ( 0.4*\tw, 0.4*\tw) rectangle ( 0.6*\tw, 0.6*\tw);
			\foreach \i in {1,...,4}
			{
				\draw[] (\i*0.2*\tw,0) -- (\i*0.2*\tw,\tw);
				\draw[] (0,\i*0.2*\tw) -- (\tw,\i*0.2*\tw);
			}
		\end{tikzpicture}
		\caption{Topology}
		\label{fig:Optimization:Topology}
    \end{subfigure} \hfill
    \caption{Types of structural optimization.}
    \label{fig:Optimization:Types}
\end{figure}

Structural optimization can be classified in three different methods: parameter, shape and topology optimization. 
Parameter optimization (\cref{fig:Optimization:Parameter}) is applied to find optimal parameters, such as the width, height or angle of an object (see e.g. \citep{Mun_2008xx,Hao_2020xx,Le_2020xx}).
Shape optimization (\cref{fig:Optimization:Shape}) is useful when parameters are not able to describe the required changes in the geometry. In the case of Computer Aided Design (CAD) representations of the geometry, control points can for example be used to adapt the surface shape (see e.g. \citep{WEEGER201926}).
Topology optimization (\cref{fig:Optimization:Topology}) allows for changing material properties and distribution to include e.g. holes into the design. Its high number of design variables offers more freedom but may limit the usability. Therefore, it is mostly applied in the initial phase to generate conceptual designs (see e.g. \citep{Lucchini_2017xx,Krenn2023xx}).

On the contrary, parameter and shape optimization are used in a later product development stage to enhance already existing designs \citep{Saitou_2005xx}. Because both parameters and shape influence the performance, a technique for coupling parameter and shape optimization is required to avoid a repetitive optimization process. Note that one could also apply shape optimization for changes in geometric parameters, as it can be seen as a more general approach than parameter optimization. However, this approach would require an excessive number of constraints to generate a practical design, e.g., to ensure rectangular permanent magnets.
One problem that arises with a coupled optimization approach is that there are many design variables. This makes the standard optimization approaches based on surrogate models or genetic algorithms computationally very expensive. Such methods have become the most widely used techniques in the last decade \citep{Duan_2013xx} but scale unfavorably with the number of design variables.

Gradient-based optimization, on the other hand, is difficult to apply with conventional mesh based optimization techniques. This is because the geometry must either be remeshed or the existing mesh must be deformed accordingly, which is especially challenging for geometry parameters that involve large deformations. This procedure has therefore been mostly restricted to shape optimization with small deformations.

This paper presents a novel method of how to combine parameter and shape optimization in an efficient way using a CAD based geometry representation. The optimization is gradient-based and uses an Isogeometric Analysis (IGA) framework with splines instead of standard tetrahedral Finite Element (FE) meshes. As a consequence, no remeshing is required since geometry changes can be performed by moving the control points of the spline curves.
The process is exemplarily demonstrated on a permanent magnet synchronous motor (PMSM) with 2D magnetostatic simulations.

The structure of the paper is as follows: First, the mathematical and physical background is explained in \cref{sec:Methodology}. Insights about the optimization algorithms are then given in \cref{sec:Optimization}. \cref{sec:Results} shows results, where the optimization is applied to a PMSM. The work in concluded in \cref{sec:Conclusions}.

\section{Methodology}
\label{sec:Methodology}

\subsection{Geometry representation}
When simulating complex designs, a suitable geometric representation is required. The use of spline-based simulation techniques has recently attracted attention as it bridges the gap between CAD programs and simulation. The flexibility of splines allows complex geometries, e.g., involving conic sections, to be accurately represented with moderate effort. This is particularly useful for rotating machines and curved surfaces. While various spline techniques exist, we restrict ourselves to the use of standard B-Splines and Non-Uniform Rational B-Splines (NURBS) \citep{Piegl_1997aa}.

\begin{figure*}[h]
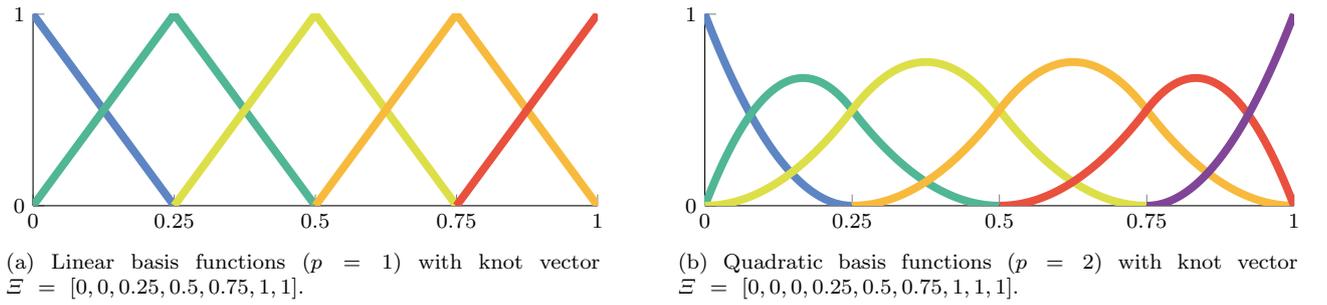

    \centering
    \begin{subfigure}{.46\textwidth}
      \centering
        \input{plots/BsplineBasis1}
        \centering
        \caption{Linear basis functions ($p=1$) with knot vector $\Xi~=~[0,0,0.25,0.5,0.75,1,1]$.}
      \label{fig:modeling:magnetocalorics:Okamura1}
    \end{subfigure} \hfill
    \begin{subfigure}{.48\textwidth}
      \centering
        \input{plots/BsplineBasis}
      \caption{Quadratic basis functions ($p=2$) with knot vector $\Xi~=~[0,0,0,0.25,0.5,0.75,1,1,1]$.}
      \label{fig:modeling:magnetocalorics:Okamura2}
    \end{subfigure}
    \centering
    \caption{Exemplary B-Spline basis functions in 1D. The knot vector indicates the support of the basis functions, i.e., a different basis function begins or ends for each entry.}
    \label{fig:Methodology:BasisFunctions}
\end{figure*}

Univariate B-Spline basis functions $\BSplineBasis_i^p$ with degree $\NurbsDegree$  are defined by the knot vector 
\begin{equation}
    \Xi =\big\{ \xi _{1} ,\ \xi _{2} ,\ ...,\ \xi _{n+p+1} \big\},
\end{equation}
with the knots $\xi_i$ which determine the support of the basis functions. Cox-de Boor's recursion formula \citep{de-Boor_1972aa}, defined as
\begin{equation}
    \BSplineBasis_{i}^{p}( \xi ) \ =\ \frac{\xi -\xi _{i}}{\xi _{i+p} -\xi _{i}} \BSplineBasis_{i}^{p-1}( \xi ) \ +\ \frac{\xi _{i+p+1} -\xi }{\xi _{i+p+1} -\xi _{i+1}} \BSplineBasis_{i+1}^{p-1}( \xi ),
\end{equation}
\begin{equation}
    \BSplineBasis_{i}^{0}( \xi ) \ =\ \begin{cases}
1 & \mathrm{if} \ \xi _{i} \leq \xi < \xi _{i+1}\\
0 & \mathrm{otherwise},
\end{cases}
\end{equation}
can then be used for the calculation of the basis functions. The number of basis functions $\NumberBasisFunctions$ depends on $\NurbsDegree$ and $\Xi$. \cref{fig:Methodology:BasisFunctions} shows exemplary basis functions of degree $p=1$ and $p=2$ in the 1D case. 

\begin{figure}
    \centering
    \input{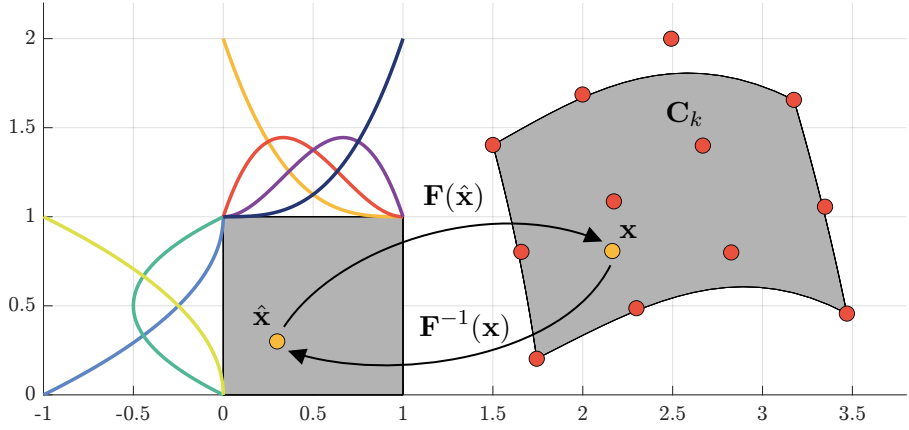}
    \caption{NURBS mapping of a two-dimensional surface from the parametric domain to the physical domain. The control points $\mathbf{\OptiCtrlPoint}_k$ in red prescribe the mapping.}
    \label{fig:Methodology:Mapping}
\end{figure}

NURBS basis functions $ \NurbsBasis_{i}^{\NurbsDegree}$ are a more general version of B-Splines and constructed via a weighting of each  $\BSplineBasis_i^p$ with an additional weight $\NurbsWeight_i$ as 
\begin{equation}
    \NurbsBasis_{i}^{\NurbsDegree}( \hat{x} ) = \frac{\NurbsWeight_{i} \BSplineBasis_{i}^{\NurbsDegree}( \hat{x} )}{\sum_{j=1}^\NumberBasisFunctions \NurbsWeight_{j} \BSplineBasis_{j}^{\NurbsDegree}( \hat{x} )}. \label{eq:NumericalMethods:Nurbs:NURBS}
\end{equation}
The rational nature of NURBS allows for an exact representation of conical sections such as circles. Multivariate B-Splines or NURBS are created in a straightforward way by using a tensor-product approach, i.e., multiplication of basis functions that are defined in each parametric direction.
To describe the mapping from the parametric space $\hat{\Omega} = (0,1)^d$ to the physical space $\Omega \subset \mathbb{R}^r$ with $r\geq d$, each geometry basis function $\hat{\GeometryBasis}$ (given as B-Spline or NURBS) is associated with a control point $\OptiCtrlPoint$ via multiplication. 
This results in the mapping
\begin{equation}
    \mathbf{F}:\hat\Omega\to\Omega, \quad
    \mathbf{x} =\mathbf{F}(\hat{\mathbf{x}}) =\sum\nolimits _{k}\hat{\GeometryBasis}_{k}(\hat{\mathbf{x}})\mathbf{\OptiCtrlPoint}_{k}. 
    \label{eq:NumericalMethods:Nurbs:3DCurve}
\end{equation}
An exemplary mapping for a two-dimensional patch is given in \cref{fig:Methodology:Mapping}. The parametric domain is shown on the left side with the basis functions defined in each dimension. The physical domain on the right side is created via the mapping \cref{eq:NumericalMethods:Nurbs:3DCurve} by multiplying the basis functions with the red control points.
Changes in the geometry are straightforwardly achieved by adjustments of the control points $\mathbf{\OptiCtrlPoint}_k$, which may lie outside the patch domain. More complex geometries -- e.g. domains which contain different materials or cannot be represented by regular transformation from the reference domain $\hat{\Omega}$ -- are constructed by a multipatch-approach, i.e., by splitting larger domains into multiple patches \citep{Buffa_2015aa}.

\subsection{Physical modeling}
\label{subsec:PhysicalModeling}
In engineering applications, the simulation of electromagnetic fields is enabled by solving Maxwell's equations. In the case of low frequency applications, such as electric motors, the magnetostatic vector potential formulation 
\begin{equation}
    \nabla\times\left(\nu\nabla\times\MagVecPot\right)= \CurrentDensity_\source+\nabla\times\left(\nu\Remanence\right)
    \label{eq:Maxwell1}
\end{equation}
is typically used, where $\MagFluxDensity = \nabla\cross\MagVecPot$ describes the magnetic flux density, $\nu = \nu\left(\norm{\MagFluxDensity}\right)$ the nonlinear reluctivity, $\CurrentDensity_\source$ the current source density and $\Remanence$ the remanence of the permanent magnets \citep{Salon_1995aa,Jackson_1998aa}. 
This formulation neglects displacement and eddy currents and is the common approximation for the case of laminated PMSMs. A further simplification is achieved by reducing \eqref{eq:Maxwell1} to the planar 2D case if a sufficiently large axial motor length is assumed. 
The computational domain is further split into the rotor $\OmegaRotor$ and the stator $\OmegaStator$ domains. This results in the Poisson problem
\begin{equation}
    \begin{cases}
    \nabla \cdotp ( \nu \nabla \MagVecPotzRt) =\nu \nabla \cdotp \RemanenceRed & \mathrm{in}\ \OmegaRotor\\
    \nabla \cdotp ( \nu \nabla \MagVecPotzSt) =-\CurrentDensityTwoD & \mathrm{in}\ \OmegaStator,
    \end{cases} \label{eq:StrongRotorStator}
\end{equation}
where only the $\MagVecPotz$-component for the magnetic vector potential is required. Note that we assume magnets in the rotor and currents in the stator. In other cases (such as e.g. reluctance machines) \cref{eq:StrongRotorStator} needs to be adapted accordingly.

The boundary conditions are indicated in 
\cref{fig:JMAGmotor}, where Dirichlet conditions on $\Gamma_\mathrm{d}$ are given by $\MagVecPotz =0$ and antiperiodic conditions are imposed on $\Gamma_\mathrm{ap}$.
For the coupling interface $\GammaAirGap$ in the air gap between rotor and stator, the  coupling condition
\begin{equation}
    \begin{cases}
    \MagVecPotzSt( \theta ) =\MagVecPotzRt( \theta - \RotAngle ) & \mathrm{on} \ \GammaAirGap \\
    H_{\theta ,\st}( \theta ) =H_{\theta ,\rt}( \theta -\RotAngle ) & \mathrm{on} \ \GammaAirGap
    \end{cases}
\end{equation}
must be fulfilled, where $\mathbf{H}=\nu\mathbf{B}$ describes the magnetic field strength. $\MagVecPotz$ and $H_\theta$ are evaluated in the respective local coordinate system fixed to rotor or stator (see  \citep{Egger_2022ab}). This corresponds to the continuity of the magnetic vector potential $\MagVecPotz$ and the azimuthal magnetic field strength $H_\theta$. From the 2D reduction follows the reduced remanence
\begin{equation}
    \RemanenceRed = \begin{pmatrix}
-B_{\mathrm{r} y}\\
B_{\mathrm{r} x}
\end{pmatrix} =B_{\mathrm{r}}\begin{pmatrix}
-\sin( \MagnetAngle )\\
\cos( \MagnetAngle )
\end{pmatrix}, \label{eq:Methodology:RemanenceDefinition}
\end{equation}
which is described by the remanence magnitude $B_\mathrm{r}$ and the direction angle $\MagnetAngle$. The source current  density in the $k$-th winding is given by
\begin{align}
\CurrentDensityTwoD^{(k)} & =\frac{\Iapp \nWind}{\Acoil}\sin\left(\polePair\RotAngle +\CurrentAngle +\frac{2\pi }{3} k\right) & k\in \{0,1,2\} \label{eq:Methodology:CurrentDefinition}
\end{align}
and zero outside the $k$-th winding. It depends on the applied current $\Iapp$, the number of windings $\nWind$, the coil cross section area $\Acoil$ and the phase number $k$.
The overall current density is then given by $\CurrentDensityTwoD=\sum\nolimits_{k}\CurrentDensityTwoD^{(k)}$ and assumed to be homogeneous within the coils. The operation of the motor is assumed to be synchronous and the rotation is defined by the electric phase offset $\CurrentAngle$ and the rotation angle $\RotAngle$ multiplied with the pole pair number $\polePair$. The time dependency in \cref{eq:Methodology:CurrentDefinition} is therefore implicitly given by $\RotAngle$.

\subsection{Numerical modeling with Isogeometric Analysis}
CAD software commonly uses splines for the representation of the geometry. The idea of IGA is to use the spline basis that is used in the CAD representation also as the basis functions for the solution. Following this approach has several intrinsic advantages, which is why the application of IGA has become increasingly popular in recent years \citep{Hughes_2005aa,Bontinck_2017ag}. The fact that conic sections can be represented exactly is especially useful in the context of rotating machines due to the circular machine geometry. Another positive aspect is that the geometry does not need to be remeshed when geometry parameters change. This has led to successful applications of IGA for simulation \citep{Bontinck_2018ac} and optimization \citep{Merkel_2021ab,Wiesheu_2023xx}.

In electrical engineering applications, the ansatz and test functions for the FE discretization are constructed with B-Splines or NURBS. Using the same functions -- here defined as $\ansatzFunction (\mathbf{x})$ -- for the ansatz and test function space yields the standard Ritz-Galerkin approach with
\begin{align}
\MagVecPotz(\mathbf{x}) \approx \solutionSymbol(\mathbf{x}) &= \sum_\ansatzFunctionIndex  \ansatzFunction_\ansatzFunctionIndex (\mathbf{x}) \solutionSymbol_\ansatzFunctionIndex & 
\testSymbol (\mathbf{x}) &= \sum_\testFunctionIndex \testFunction_\testFunctionIndex (\mathbf{x})\testSymbol_\testFunctionIndex. \label{eq:AnsatzTestFunctions}
\end{align}

After derivation of the weak formulation of \eqref{eq:StrongRotorStator} and using \eqref{eq:AnsatzTestFunctions} for the discretization, we obtain the  coupled matrix system

\begin{equation}
    \underbrace{\begin{pmatrix}
    \stiffnessMatrix_\rt & \mathbf{0} & -\couplingMatrix_\rt \\
    \mathbf{0} & \stiffnessMatrix_\st & \couplingMatrix_\st\RotMat \\
    -\couplingMatrix_\rt^\top & \RotMat^\top\couplingMatrix_\st^\top & \mathbf{0}
    \end{pmatrix}}_{=:\stiffnessMatrix}\underbrace{\begin{pmatrix}
    \solutionVector_\rt \\
    \solutionVector_\st \\
    \mortarVector
    \end{pmatrix}}_{=:\solutionVector} =\underbrace{\begin{pmatrix}
    \rhsVector_\rt \\
    \rhsVector_\st \\
    \mathbf{0}
    \end{pmatrix}}_{=:\rhsVector} \label{eq:Methodology:CoupledMatrixSystem}
\end{equation}
with the matrix entries
\begin{align}
K_{ij} & =\int _{\Omega } \nu \nabla N_{i}(\mathbf{x}) \cdotp \nabla N_{j}(\mathbf{x})\der \Omega \label{eq:methodology:Krt} \\ 
b_{i,\rt} & =\int _{\Omega } \nu \nabla N_{i}(\mathbf{x}) \cdotp \mathbf{B}_{\mathrm{r}}^{\bot }\der \Omega \label{eq:methodology:brt} \\ 
b_{i,\st} & =\int _{\Omega } J_{z,\source} N_{i}(\mathbf{x})\der \Omega \label{eq:Methodology:bst}
\end{align} 
for both the rotor and stator domains. The entries for the coupling matrices $\couplingMatrix_\rt$ and $\couplingMatrix_\st$ arise from the usage of harmonic (Fourier) functions for the mortaring on the interface of rotor and stator and are given as
\begin{align}
g_{i,2n-1} & =\int _{\Gamma _{\ag}} N_{i}\sin( n\theta )\der \Gamma  & n\geq 1\\
g_{i,2n} & =\int _{\Gamma _{\ag}} N_{i}\cos( n\theta )\der \Gamma  & n\geq 0
\end{align}
where $\GammaAirGap$ is evaluated for the rotor and stator side, respectively. For different rotation angles $\RotAngle \in \{\RotAngle_1,\RotAngle_2,...,\RotAngle_N \}$, the stator coupling matrix $\couplingMatrix_\st$ is multiplied with a rotation matrix $\RotMat$, which contains blockdiagonal sine and cosine entries, such that the coupling matrix does not need to be reassembled for each rotation angle. This makes evaluations for multiple rotation angles very efficient, as only $\RotMat$ needs to be updated. Further details are found in \citep{Egger_2022ab}. Due to the symmetry of the motor, only one quarter is simulated with antiperiodic boundary conditions on the interface $\Gamma_\mathrm{ap}$, see \cref{fig:JMAGmotor}. Because of this, only the harmonics $n\in \{2,6,10,...\}$ are non-zero.

Overall, \eqref{eq:Methodology:CoupledMatrixSystem} can be written as
\begin{equation}
\stateVector(\solutionVector) = \stiffnessMatrix(\solutionVector)\solutionVector - \rhsVector = \mathbf{0}, \label{eq:Methodology:CoupledMatrixSystemred}
\end{equation}
where $\stateVector$ is introduced as variable for the state equation. Note that $\stiffnessMatrix$ depends on the solution of the magnetic flux density due to the nonlinear material law. We solve \eqref{eq:Methodology:CoupledMatrixSystemred} therefore iteratively with the Newton-Raphson scheme, i.e. by linearizing \eqref{eq:Methodology:CoupledMatrixSystemred} using the derivatives given in \ref{app:NewtonJacobian}.
The torque is calculated with
\begin{equation}
    T_\RotAngle(\solutionVector) = -L\solutionVector_\st^\top \couplingMatrix_\st \RotMat^\prime \mortarVector, \label{eq:Methogology:Torque}
\end{equation}
which is derived from energy balance considerations \citep{Egger_2022ab}. 
As the computations are carried out in 2D, the axial length $L$ of the motor is explicitly considered in \cref{eq:Methogology:Torque}. The quantity $\RotMat^\prime$ describes the derivative of $\RotMat$ with respect to $\RotAngle$.

\section{Optimization}
\label{sec:Optimization}
The simulation and optimization of electric motors is a complex task with high industrial relevance. 
For a comprehensive study of all relevant physical effects, one has to study among others electromagnetic fields, mechanical stresses and thermal effects \citep{Rosu_2017aa}. 
For demonstration purposes of our optimization method, we restrict ourselves here to the magnetic aspects.

\begin{figure*}[t]
    \centering
    \input{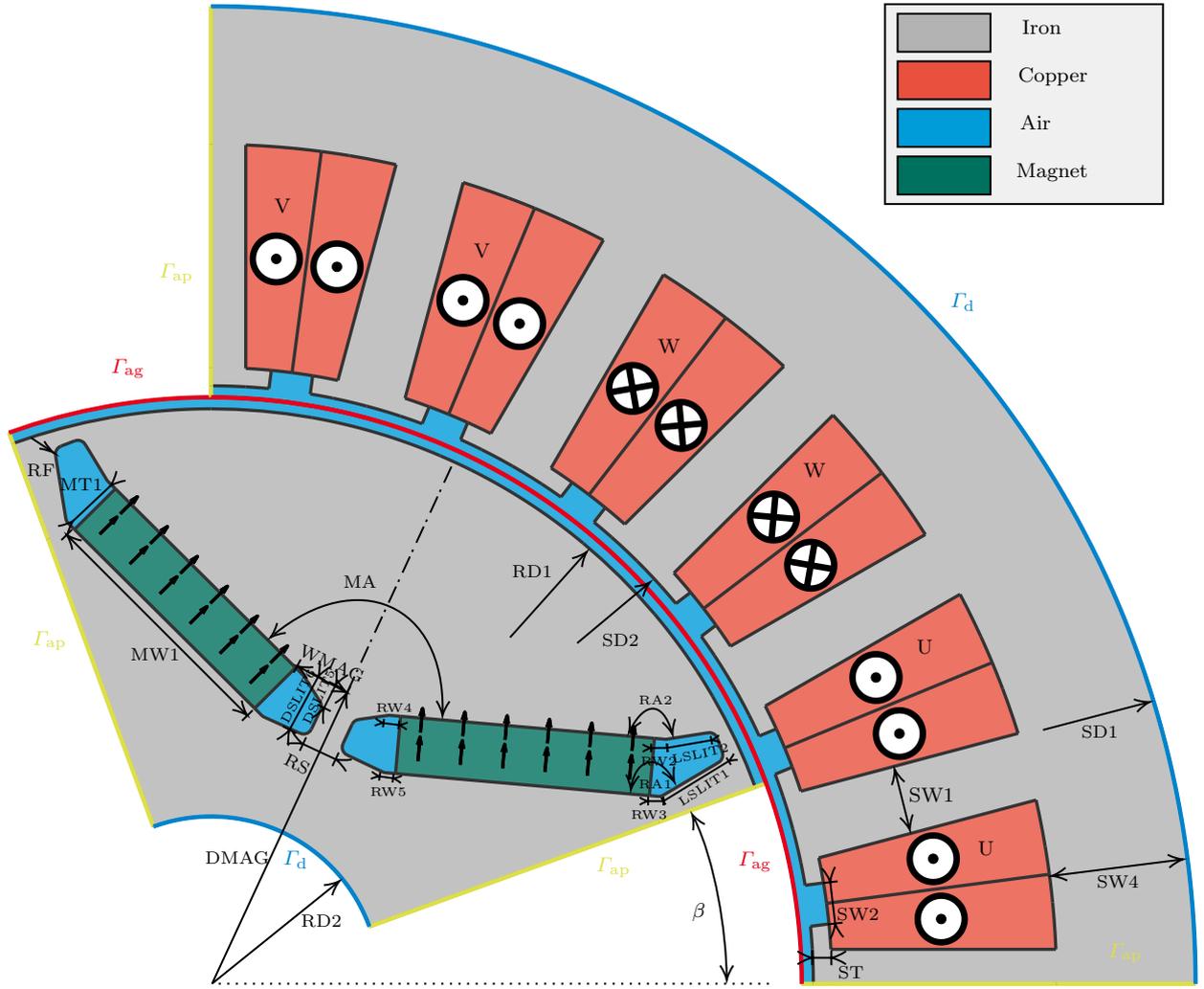}
    \caption{PMSM geometry, rebuilt from JMAG. Geometry parameters, material definitions and boundary conditions are provided.}
    \label{fig:JMAGmotor}
\end{figure*}

This paper uses a benchmark PMSM geometry from the commercial program \textit{JMAG}\textsuperscript{\textregistered}. The 2D geometry is rebuilt in \textit{MATLAB}\textsuperscript{\textregistered} with a spline patch representation using the open source software \textit{GeoPDEs} \citep{Vazquez_2016aa} and the \textit{NURBS} toolbox \citep{NURBS_2021aa}. \cref{fig:JMAGmotor} shows a detailed picture of the motor including geometry parameters, material distribution and boundaries. The used geometric parameters are given in \cref{app:tab:Parameters}. The green patches indicate permanent magnets, where the magnetization direction is given by the arrows. Gray and blue patches represent iron and air, respectively. The windings are depicted by red patches with the winding scheme as indicated. 

\subsection{Problem definition}
When optimizing electric machines there are several relevant quantities to consider. 
For commercial applications, the minimization of production cost is one key interest. This is achieved primarily by reducing material usage. In the case of rare-earth materials, this becomes even more important due to environmental reasons \citep{Binnemans_2013xx}. Another important feature is that the motor fulfills the specified performance requirements, where we aim to maintain the mean torque $\MeanTorque$. Furthermore, the qualitative profile of the torque is important, as torque ripple causes unwanted vibrations and noise. Therefore, several works address this problem \citep{Sanchez_2022xx}. Generating a uniform torque profile has been achieved in several ways, e.g. by restricting the interval of maximum and minimum torque values \citep{Brun_2023aa,Zhentao_2017xx},  minimization of the torque standard deviation $\StdTorque$ \citep{Li_2015xx}, as well as by considering the total harmonic distortion of the torque \citep{Bontinck_2017ai}, the electromotive force (EMF) \citep{Merkel_2021ab} or magnetic flux density \citep{Choi_2012xx,Xie2019xx}. So far, this is obtained by either parameter or shape optimization. There is little research to combine different types of optimization \citep{Kuci_2019aa,Guo_2022aa,Weeger2022xx}. To increase the design freedom, a combined approach is proposed in the following.

We propose to use a multi-objective optimization approach with a weighting of the different objectives. This makes the optimization customizable to the user's needs and reduces the computational cost compared to classical multi-objective optimization, but is limited to cases with convex Pareto fronts \citep{Ehrgott_2005aa}.  
We include the area of the permanent magnets $A_\mathrm{Magnet}$, the standard deviation of the torque $\StdTorque$ and a regularization term $S$ for the surface smoothness as objective functions to be minimized. 
Calculating $\MeanTorque$ and $\StdTorque$ for a discrete set of rotation angles $\RotAngle$ is achieved with
\begin{align}
    \MeanTorque &=\frac{1}{N}\sum\nolimits _{\RotAngle } T_{\RotAngle } & \StdTorque &=\sqrt{\frac{1}{N}\sum\nolimits _{\RotAngle }( T_{\RotAngle } -\MeanTorque)^{2}}.
    \label{eq:Optimization:Torques}
\end{align}
Reducing $A_\mathrm{Magnet}$ is desired because permanent magnets are the most expensive parts in PMSMs due to the necessary rare earth elements. Minimizing $A_\mathrm{Magnet}$  will therefore yield a low cost design and reduce the negative ecological impact. We choose to reduce the ripples by minimizing $\StdTorque$, which is continuously differentiable and efficient to compute in the current framework as will be explained later.
To ensure the required performance, $\MeanTorque$ is bounded from below by the minimum acceptable torque $T_\mathrm{Target}$.

The geometry is initially built with parameters  $\OptiParameter$ (length, width, radii, ...) which give rise to control points $\OptiCtrlPoint_\mathrm{init}$ that are then adjusted with radial offsets $\mathrm{\Delta}\OptiCtrlPoint$ for the rotor surface, i.e.,
\begin{equation}
    \OptiCtrlPoint(\OptiParameter , \mathrm{\Delta}\OptiCtrlPoint) = \OptiCtrlPoint_\mathrm{init}(\OptiParameter) + \mathrm{\Delta}\OptiCtrlPoint.
\end{equation}
This definition allows for convenient surface modifications and calculation of gradients as $\mathrm{d}\OptiCtrlPoint/\mathrm{d}\mathrm{\Delta}\OptiCtrlPoint=1$. The derivatives with respect to $\OptiCtrlPoint$ can therefore also be used for $\mathrm{\Delta}\OptiCtrlPoint$. 
The regularization term $S$ is calculated with 
\begin{equation}
    S=\sum\nolimits _{i}\frac{( \mathrm{\Delta} C_{i+1} -\mathrm{\Delta} C_i)^{2}}{\theta _{i+1} -\theta _{i}}
\end{equation}
as the squared difference of consecutive control point offsets $\mathrm{\Delta} C$, divided by the difference in their angular position~$\theta$. This corresponds to a regularization of the rotor surface by a squared discrete $H^1$-norm.
As design variables, the parameters  $\OptiParameter$ and the control point offsets $\mathrm{\Delta}\OptiCtrlPoint$ are chosen. 
The magnetic field solution $\solutionVector$ depends therefore on the control points $\OptiCtrlPoint$ and parameters $\OptiParameter$, i.e.,
\begin{equation}
    \solutionVector =\solutionVector\bigl(\OptiCtrlPoint(\OptiParameter , \mathrm{\Delta}\OptiCtrlPoint) , \OptiParameter\bigr),
\end{equation}
where $\OptiCtrlPoint$ again depends on the geometry parameters in $\OptiParameter$ and subsequent shape modifications performed with $\mathrm{\Delta}\OptiCtrlPoint$. The overall problem is then defined as
\begin{equation}
    \begin{array}{l}
        \min \fopt(\optiVector,\solutionVector) = w_1\, \Amagnet + w_2\, \StdTorque + w_3 S \\[2.5pt]
        \mathrm{s.t.} \\ 
        \begin{cases}
        \hspace{0.5cm} \stateVector(\optiVector,\solutionVector) = \mathbf{0} & \mathrm{State~equation}\\
        \hspace{0.5cm} \MeanTorque \geq T_{\mathrm{Target}} & \mathrm{Fulfill~target~torque} \\
        \hspace{0.5cm} \mathbf{0} \leq \optiVector \leq \mathbf{1} & \mathrm{Design~limits} \\
        \hspace{0.5cm} \mathbf{g}(\optiVector) \leq \mathbf{0} & \mathrm{Geometric~feasibility} 
        
        \end{cases}
    \end{array}\label{eq:Optimization:Problem}
\end{equation}
where in the design vector $\optiVector$ the values for $\OptiParameter$ and $\mathrm{\Delta}\OptiCtrlPoint$ are concatenated with an additional scaling given by the min and max values such that $\mathbf{0} \leq \optiVector \leq \mathbf{1}$.
This makes \eqref{eq:Optimization:Problem} a PDE-constrained optimization problem, where it is also necessary to solve \eqref{eq:Methodology:CoupledMatrixSystemred}. Note that the state equation $\stateVector(\optiVector,\solutionVector)$ has now also the additional dependency on~$\optiVector$.
Additionally, in \eqref{eq:Optimization:Problem}, geometric constraints $\mathbf{g}(\optiVector)$  have been added in order to ensure a feasible design. The inequality there has to be read component-wise.

The quantities $\MeanTorque$ and $\StdTorque$ are computed by evaluation of \eqref{eq:Methogology:Torque} at a discrete set of rotation angles.
The weights $w_1$, $w_2$ and $w_3$ are introduced in order to stabilize the problem by weighting the optimization objectives such that the contributions are of comparable magnitude.

\subsection{Sensitivity analysis}
When optimizing computationally expensive problems with a large number of optimization variables, employing a gradient-based algorithm is attractive to avoid prohibitive computation times.  Numerical approximation of the gradients, e.g. by finite differences, scale with the number of design variables and is often unstable. For an efficient optimization process, it is therefore necessary to resort to analytical derivations for the gradients whenever possible. As we combine parameter and shape optimization, the design space $\optiVector$ is comparatively large.

The most challenging derivatives to obtain are the ones of the motor's torque $\Torque_{\RotAngle}$ at a given rotation angle $\RotAngle$. For that, first the derivatives of $\Torque_{\RotAngle}$ are derived in \cref{subsubsec:ShapeSensitivities} for the shape and in \cref{subsubsec:ParameterSensitivities} for the parameters. Finally, the sensitivities of $\MeanTorque$ and $\StdTorque$ are then derived in \cref{subsubsec:PostprocessingTorqueDer}. 

\subsubsection{Shape optimization sensitivities}
\label{subsubsec:ShapeSensitivities}
This paper uses a \textit{discretize-then-optimize} approach. This means that the optimization gradients are calculated after the discretization of the finite dimensional problem.  
Consequently, the derivatives are built with respect to the control points $\OptiCtrlPoint$, which define the geometry via the mapping \eqref{eq:NumericalMethods:Nurbs:3DCurve}. This is in contrast to the \textit{optimize-then-discretize} approach, where the gradients are calculated in the continuous setting and discretized afterwards, which has been also applied in the IGA context \citep{Merkel_2021ab,Ziegler_2023xx}.
Note that both approaches have been shown to be identical for simple settings \citep{Fuseder_2015aa}.
This is a feature of IGA being able to represent the geometry exactly.

In a  discrete setting,  the derivatives of a function $f$ with respect to the design $\optiVector$ can be accomplished in a general way by solving
\begin{equation}
    \frac{\der f(\optiVector ,\solutionVector)}{\der\optiVector} =\frac{\partial f(\optiVector ,\solutionVector)}{\partial \optiVector} +\left(\frac{\partial \stateVector(\optiVector ,\solutionVector)}{\partial \optiVector}\right)^{\top }\mathbf{\adjointVector }
    \label{eq:generalderivative}
\end{equation}
as derived from \citep{Hinze_2009aa}, which is called the adjoint method. The adjoint solution $\adjointVector$ is obtained from solving
\begin{equation}
    \frac{\partial f(\optiVector ,\solutionVector)}{\partial \solutionVector} +\left(\frac{\partial \stateVector(\optiVector ,\solutionVector)}{\partial \solutionVector}\right)^{\top }\mathbf{\adjointVector } =\mathbf{0}.
\end{equation}
Applying this procedure to the derivative of $\Torque_\RotAngle$ with respect to $\ControlPointSymbol_{kd}$ ($k$'th control point in dimension $d$) yields
\begin{equation}
    \frac{\der \Torque_\RotAngle}{\der \ControlPointSymbol_{kd}} =\frac{\partial \Torque_\RotAngle}{\partial \ControlPointSymbol_{kd}} +\sum\nolimits _{i} \adjointSymbol_{i}\left(\sum\nolimits _{j}\frac{\der \stiffnessSymbol_{ij}}{\der \ControlPointSymbol_{kd}} u_{j} -\frac{\der \rhsSymbol_{i}}{\der \ControlPointSymbol_{kd}}\right), \label{eq:Optimization:FullDerivative}
\end{equation}
where the adjoint variable $\adjointVector$ is computed by solving
\begin{equation}
    \left(\frac{\der\mathbf{\stiffnessMatrix}}{\der\mathbf{\solutionVector}}\mathbf{\solutionVector} +\mathbf{\stiffnessMatrix}\right)^{\top }\mathbf{\adjointVector }
    =-\frac{\partial \Torque_\RotAngle}{\partial \solutionVector}. \label{eq:Optimization:AdjointMethod}
\end{equation}
The left hand side of \eqref{eq:Optimization:AdjointMethod} is given by the Jacobian of $\stateVector(\solutionVector)$ with respect to $\solutionVector$. The product is understood as
\begin{equation}
    \left(\frac{\der\stiffnessMatrix}{\der\solutionVector}\solutionVector \right)_{ij} = \sum_k\frac{\der\stiffnessSymbol_{ik}}{\der\solutionSymbol_j}\solutionSymbol_k.
\end{equation}

This means that five derivatives must be obtained for an analytical expression of \eqref{eq:Optimization:FullDerivative}. The first expression $\frac{\partial \Torque_\RotAngle}{\partial \ControlPointSymbol_{kd}}$ is zero, as $\fopt$  does not depend explicitly  on $\ControlPointSymbol_{kd}$. The derivative $\frac{\partial \Torque_\RotAngle}{\partial \solutionVector}$ is obtained by differentiating its analytical expression \eqref{eq:Methogology:Torque} and leads to
\begin{equation}
    \frac{\partial T_\RotAngle}{\partial \solutionVector}=-L\left(\begin{array}{c}
    \mathbf{0}\\
    \couplingMatrix_\st \RotMat^{\prime}\mortarVector\\
    \RotMat^{\prime\top}\couplingMatrix^\top_\st\solutionVector_\st
    \end{array}\right),
\end{equation}
where $\RotMat$, $\mortarVector$ and $\solutionVector_\st$ depend on the rotation angle $\RotAngle$. 

The Jacobian of $\stateVector(\solutionVector)$ needs to be evaluated regardlessly during the Newton-Raphson scheme for solving \eqref{eq:Methodology:CoupledMatrixSystemred} and is given in \ref{app:NewtonJacobian}. 

Two additional derivatives  $\frac{\der \stiffnessSymbol_{ij}}{\der \ControlPointSymbol_{kd}}$ and $\frac{\der \rhsSymbol_{i}}{\der \ControlPointSymbol_{kd}}$ need to be determined.  They quantify how the stiffness matrix and right hand side of \eqref{eq:Methodology:CoupledMatrixSystem} change, when the $k$-th control point is moved in direction $d$ (where $d=x$ or $d=y$).
To avoid having to compute gradients of the coupling matrices $\couplingMatrix_\rt$ and $\couplingMatrix_\st$, we keep the discretization at $\GammaAirGap$ unchanged, since the optimization only changes the interior domain and not the boundary discretization. Using the transformation rule $\nabla N_{i} =\mathbf{J}_{F}^{-\top }\hat{\nabla }\hat{N}_{i}$ for functions in $H^1$ from \citep{Monk_2003aa} on \eqref{eq:methodology:Krt} and \eqref{eq:methodology:brt} yields
\begin{align}
    K_{ij,\rt} &=\int_{\hat{\Omega }} \nu(B) \left(\mathbf{J}_{F}^{-\top }\hat{\nabla }\hat{N}_{i}\right) \cdotp \left(\mathbf{J}_{F}^{-\top }\hat{\nabla }\hat{N}_{j}\right) |\mathbf{J}_{F} |\der\hat{\Omega } \label{eq:Optimization:Krt} \\
    b_{i,\rt} & = \int _{\hat{\Omega }} \nu\left(\mathbf{J}_{F}^{-\top }\hat{\nabla }\hat{N}_{i}\right) \cdotp \mathbf{B}_{\mathrm{r}}^{\bot } |\mathbf{J}_{F} |\der\hat{\Omega }, \label{eq:Optimization:brt}
\end{align}
which shows their dependency on the Jacobian matrix $\Jac$. The latter contains the mapping information and therefore depends on the control points. $\Jac$ is derived from the mapping \eqref{eq:NumericalMethods:Nurbs:3DCurve} and calculated with 
\begin{equation}
    \Jac = \begin{pmatrix}
    \sum\nolimits _{k}\left(\hat{\nabla }\hat{\GeometryBasis}_{k}\right)^{\top }\ControlPointSymbol_{kx}\\
    \sum\nolimits _{k}\left(\hat{\nabla }\hat{\GeometryBasis}_{k}\right)^{\top }\ControlPointSymbol_{ky}
    \end{pmatrix}\label{eq:Optimization:JacobianMapping}
\end{equation}
in the 2D case. It can be seen that \eqref{eq:Optimization:JacobianMapping} directly depends on the control points via multiplication. Differentiation of \eqref{eq:Optimization:Krt} and \eqref{eq:Optimization:brt} is achieved by using the product rule on all components that contain $\Jac$. These analytic expressions are shown in more detail in \ref{app:sensitivityStiffnessMat}. As a result, we obtain the derivatives $\frac{\der \stiffnessSymbol_{ij}}{\der \ControlPointSymbol_{kd}}$ and $\frac{\der \rhsSymbol_{i}}{\der \ControlPointSymbol_{kd}}$. 
One feature of these derivatives is that they have a high degree of sparsity due to the tensor product structure. 
The computation is also highly parallelizable over multiple patches, patch elements and integral contributions.

\subsubsection{Parameter optimization sensitivities}
\label{subsubsec:ParameterSensitivities}

Calculating the sensitivities of the geometry parameters is achieved in a semi-analytical way making use of the already calculated control point derivatives. The process is based on the idea that changing one geometry parameter corresponds to changing several control points at once. So if one knows the contribution $\frac{\der \ControlPointSymbol_{kd}}{\der \ParameterSymbol_{l}}$, i.e., how the control point $\ControlPointSymbol_{kd}$ changes if a geometry parameter $\ParameterSymbol_l$ is changed, the parameter sensitivities can be directly calculated from this. Mathematically speaking, we apply the chain rule and the parameter derivatives can be calculated with
\begin{equation}
    \frac{\der \Torque}{\der \ParameterSymbol_{l}} =\underbrace{\sum _{k,d}\frac{\der \Torque}{\der \ControlPointSymbol_{kd}}\frac{\der \ControlPointSymbol_{kd}}{\der \ParameterSymbol_{l}}}_{D_{\mathrm{geo}}} +\underbrace{\sum _{k}\frac{\partial \Torque}{\partial \solutionSymbol_{k}}\frac{\partial \solutionSymbol_{k}}{\partial \ParameterSymbol_{l}}}_{D_{\mathrm{phys}}}. \label{eq:ChainRule}
\end{equation}
Any contribution that involves geometric changes that can be represented by modifications of the control points are represented in the first expression of \eqref{eq:ChainRule} by $D_{\mathrm{geo}}$. As $\frac{\der \Torque_\RotAngle}{\der \ControlPointSymbol_{kd}}$ has already been derived in the previous section, only $\frac{\der \ControlPointSymbol_{kd}}{\der \ParameterSymbol_{l}}$ must be additionally computed. This would be possible in an analytic way, but very tedious due to the different patch types (here B-Splines and NURBS) and the various ways of how parameters can define patches (e.g. lengths, radii or angles). Instead, a numerical approach is used, where the motor geometry is rebuilt for each parameter with perturbed values. A schematic explanation is given in \cref{fig:Optimization:DPDC} for a single patch geometry with the initial parameter $P_0$ and an exemplary control point $C_0$. The derivative is then obtained by a perturbed geometry with parameter $P_1$ and the corresponding $C_1$ using e.g. forward differences, i.e., $(C_1-C_0)/(P_1-P_0)$.
\begin{figure*}[h]
    \centering
    \input{plots/DCDP.tikz}
    \caption{Visualization of the numerical evaluation of control point changes depending on a parameter used in the evaluation of \eqref{eq:ChainRule}.}
    \label{fig:Optimization:DPDC}
\end{figure*}

This has the advantage that it is generally applicable and simple to implement. Due to the drastically reduced number of patches (compared to conventional FE), this method is also fast, e.g. in our case it is only a matter of milliseconds per parameter. Numerical stability is also not a problem, since the admissible parameter ranges are known and the step size can be chosen accordingly. Finally, this method allows for a very favorable scaling in the number of geometry parameters, as the most time-consuming step is building the derivatives from the previous section.

Note that this approach neglects the influence of the weights $\NurbsWeight$ in case of a NURBS geometry representation from \eqref{eq:NumericalMethods:Nurbs:NURBS}. One can either incorporate the influence analogous to \eqref{eq:ChainRule} using the analytic gradients with respect to $\NurbsWeight$ (see e.g. \citep{Qian_2010xx}) or enforce the weights to remain unchanged.
Another practical approach -- used in this work -- is to disregard the influence of the NURBS weights. This is possible as the influence is negligible compared to the one of the control points, such that there is no significant effect on the optimization.

The second contribution $D_{\mathrm{phys}}$ in \eqref{eq:ChainRule} represents physical quantities that cannot be expressed via movements of control points. In the case of an electric motor, two such quantities are the magnetization direction $\MagnetAngle$ of the permanent magnets (given in \eqref{eq:Methodology:RemanenceDefinition}) and the electric operating angle of the three-phase current $\CurrentAngle$ (given in \eqref{eq:Methodology:CurrentDefinition}). Their contribution is incorporated by the derivatives of \eqref{eq:methodology:brt} and \eqref{eq:Methodology:bst} with respect to $\MagnetAngle$ and $\CurrentAngle$ -- which are simple sine and cosine derivatives -- respectively (see \ref{app:sensitivityStiffnessMat}). $D_{\mathrm{phys}}$ is then calculated with  \eqref{eq:generalderivative} using the known adjoint solution $\adjointVector$. This allows the operating angle to be conveniently determined during the optimization.

\subsubsection{Calculation of torque derivatives}
\label{subsubsec:PostprocessingTorqueDer}
Once the derivatives of the torque for different rotation angles $\RotAngle$ are known, the mean and standard deviation are calculated with  \cref{eq:Optimization:Torques}. The derivatives of the mean torque $\MeanTorque$ and torque ripple $\StdTorque$ with respect to any parameter $\ParameterSymbol_l$ are
\begin{align}
    \frac{\der\MeanTorque}{\der\ParameterSymbol_l} &=\frac{1}{N}\sum\nolimits _{\RotAngle }\frac{\der T_{\RotAngle }}{\der\ParameterSymbol_l} \\
    \frac{\der\StdTorque}{\der \ParameterSymbol_l} &=\frac{1}{\StdTorque}\left(\frac{1}{N}\sum\nolimits _{\RotAngle } \Torque_{\beta }\frac{\der \Torque_{\RotAngle }}{\der \ParameterSymbol_l} -\MeanTorque\frac{\der\MeanTorque}{\der \ParameterSymbol_l}\right)
\end{align}
and carried out analogously for the control points. The same approach is found in the literature \citep{Li_2015xx}.

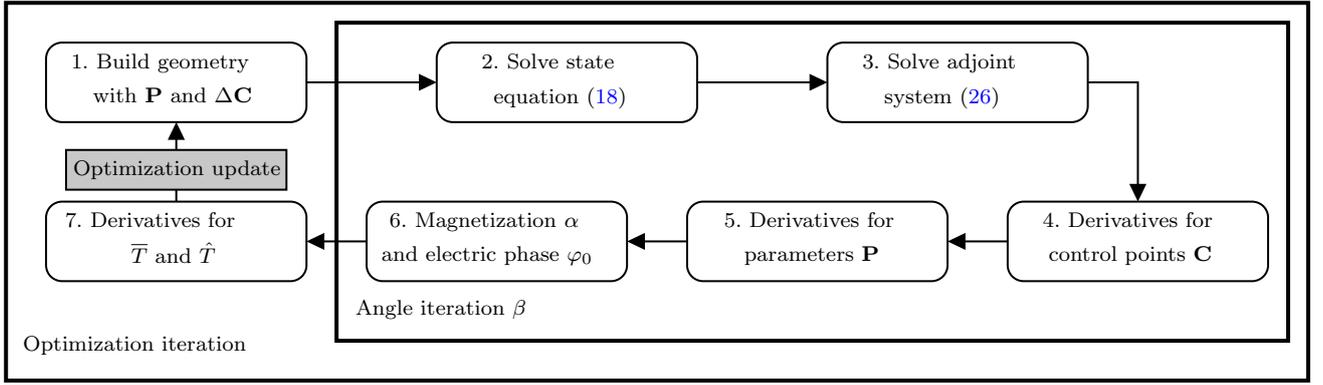
\begin{figure*}[t]
    \centering
    \tikzset{every picture/.style={line width=0.75pt}} 

\begin{tikzpicture}[x=0.75pt,y=0.75pt,yscale=-1,xscale=1]

\draw  [color={rgb, 255:red, 0; green, 0; blue, 0 }  ,draw opacity=1 ][line width=1.5]  (0,70) -- (650,70) -- (650,260) -- (0,260) -- cycle ;
\draw  [color={rgb, 255:red, 0; green, 0; blue, 0 }  ,draw opacity=1 ][line width=1.5]  (165,80) -- (640,80) -- (640,240) -- (165,240) -- cycle ;
\draw  [fill={rgb, 255:red, 255; green, 255; blue, 255 }  ,fill opacity=1 ] (20,98) .. controls (20,93.58) and (23.58,90) .. (28,90) -- (142,90) .. controls (146.42,90) and (150,93.58) .. (150,98) -- (150,122) .. controls (150,126.42) and (146.42,130) .. (142,130) -- (28,130) .. controls (23.58,130) and (20,126.42) .. (20,122) -- cycle ;
\draw    (85,210) -- (85,133) ;
\draw [shift={(85,130)}, rotate = 90] [fill={rgb, 255:red, 0; green, 0; blue, 0 }  ][line width=0.08]  [draw opacity=0] (8.93,-4.29) -- (0,0) -- (8.93,4.29) -- cycle    ;
\draw    (180,190) -- (153,190) ;
\draw [shift={(150,190)}, rotate = 360] [fill={rgb, 255:red, 0; green, 0; blue, 0 }  ][line width=0.08]  [draw opacity=0] (8.93,-4.29) -- (0,0) -- (8.93,4.29) -- cycle    ;
\draw    (380,200) -- (353,200) ;
\draw [shift={(350,200)}, rotate = 360] [fill={rgb, 255:red, 0; green, 0; blue, 0 }  ][line width=0.08]  [draw opacity=0] (8.93,-4.29) -- (0,0) -- (8.93,4.29) -- cycle    ;
\draw  [fill={rgb, 255:red, 255; green, 255; blue, 255 }  ,fill opacity=1 ] (410,98) .. controls (410,93.58) and (413.58,90) .. (418,90) -- (532,90) .. controls (536.42,90) and (540,93.58) .. (540,98) -- (540,122) .. controls (540,126.42) and (536.42,130) .. (532,130) -- (418,130) .. controls (413.58,130) and (410,126.42) .. (410,122) -- cycle ;
\draw    (540,110) -- (565,110) -- (565,167) ;
\draw [shift={(565,170)}, rotate = 270] [fill={rgb, 255:red, 0; green, 0; blue, 0 }  ][line width=0.08]  [draw opacity=0] (8.93,-4.29) -- (0,0) -- (8.93,4.29) -- cycle    ;
\draw  [fill={rgb, 255:red, 255; green, 255; blue, 255 }  ,fill opacity=1 ] (20,178) .. controls (20,173.58) and (23.58,170) .. (28,170) -- (142,170) .. controls (146.42,170) and (150,173.58) .. (150,178) -- (150,202) .. controls (150,206.42) and (146.42,210) .. (142,210) -- (28,210) .. controls (23.58,210) and (20,206.42) .. (20,202) -- cycle ;
\draw  [fill={rgb, 255:red, 255; green, 255; blue, 255 }  ,fill opacity=1 ] (180,178) .. controls (180,173.58) and (183.58,170) .. (188,170) -- (302,170) .. controls (306.42,170) and (310,173.58) .. (310,178) -- (310,202) .. controls (310,206.42) and (306.42,210) .. (302,210) -- (188,210) .. controls (183.58,210) and (180,206.42) .. (180,202) -- cycle ;
\draw  [fill={rgb, 255:red, 255; green, 255; blue, 255 }  ,fill opacity=1 ] (340,178) .. controls (340,173.58) and (343.58,170) .. (348,170) -- (462,170) .. controls (466.42,170) and (470,173.58) .. (470,178) -- (470,202) .. controls (470,206.42) and (466.42,210) .. (462,210) -- (348,210) .. controls (343.58,210) and (340,206.42) .. (340,202) -- cycle ;
\draw  [fill={rgb, 255:red, 255; green, 255; blue, 255 }  ,fill opacity=1 ] (500,178) .. controls (500,173.58) and (503.58,170) .. (508,170) -- (622,170) .. controls (626.42,170) and (630,173.58) .. (630,178) -- (630,202) .. controls (630,206.42) and (626.42,210) .. (622,210) -- (508,210) .. controls (503.58,210) and (500,206.42) .. (500,202) -- cycle ;
\draw    (340,190) -- (313,190) ;
\draw [shift={(310,190)}, rotate = 360] [fill={rgb, 255:red, 0; green, 0; blue, 0 }  ][line width=0.08]  [draw opacity=0] (8.93,-4.29) -- (0,0) -- (8.93,4.29) -- cycle    ;
\draw    (500,190) -- (473,190) ;
\draw [shift={(470,190)}, rotate = 360] [fill={rgb, 255:red, 0; green, 0; blue, 0 }  ][line width=0.08]  [draw opacity=0] (8.93,-4.29) -- (0,0) -- (8.93,4.29) -- cycle    ;
\draw  [fill={rgb, 255:red, 200; green, 200; blue, 200 }  ,fill opacity=1 ] (30,144) -- (140,144) -- (140,164) -- (30,164) -- cycle ;
\draw    (345,110) -- (407,110) ;
\draw [shift={(410,110)}, rotate = 180] [fill={rgb, 255:red, 0; green, 0; blue, 0 }  ][line width=0.08]  [draw opacity=0] (8.93,-4.29) -- (0,0) -- (8.93,4.29) -- cycle    ;
\draw  [fill={rgb, 255:red, 255; green, 255; blue, 255 }  ,fill opacity=1 ] (215,98) .. controls (215,93.58) and (218.58,90) .. (223,90) -- (337,90) .. controls (341.42,90) and (345,93.58) .. (345,98) -- (345,122) .. controls (345,126.42) and (341.42,130) .. (337,130) -- (223,130) .. controls (218.58,130) and (215,126.42) .. (215,122) -- cycle ;
\draw    (150,110) -- (212,110) ;
\draw [shift={(215,110)}, rotate = 180] [fill={rgb, 255:red, 0; green, 0; blue, 0 }  ][line width=0.08]  [draw opacity=0] (8.93,-4.29) -- (0,0) -- (8.93,4.29) -- cycle    ;

\draw (236,94) node [anchor=north west][inner sep=0.75pt]   [align=left] {2. Solve state};
\draw (242,111) node [anchor=north west][inner sep=0.75pt]   [align=left] {equation \eqref{eq:Methodology:CoupledMatrixSystemred}};
\draw (31,94) node [anchor=north west][inner sep=0.75pt]   [align=left] {1. Build geometry};
\draw (42,111) node [anchor=north west][inner sep=0.75pt]   [align=left] {with $\OptiParameter$ and $\mathrm{\Delta}\OptiCtrlPoint$ };
\draw (61,189) node [anchor=north west][inner sep=0.75pt]   [align=left] {$\MeanTorque$ and $\StdTorque$};
\draw (426,94) node [anchor=north west][inner sep=0.75pt]   [align=left] {3. Solve adjoint };
\draw (28,174) node [anchor=north west][inner sep=0.75pt]   [align=left] {7. Derivatives for};
\draw (357,174) node [anchor=north west][inner sep=0.75pt]   [align=left] {5. Derivatives for};
\draw (190,174) node [anchor=north west][inner sep=0.75pt]   [align=left] {6. Magnetization $\MagnetAngle$};
\draw (7,236) node [anchor=north west][inner sep=0.75pt]   [align=left] {Optimization iteration};
\draw (173,218) node [anchor=north west][inner sep=0.75pt]   [align=left] {Angle iteration $\RotAngle$};
\draw (516,174) node [anchor=north west][inner sep=0.75pt]   [align=left] {4. Derivatives for };
\draw (437,111) node [anchor=north west][inner sep=0.75pt]   [align=left] {system $\eqref{eq:Optimization:AdjointMethod}$};
\draw (186,191) node [anchor=north west][inner sep=0.75pt]   [align=left] {and electric phase $\CurrentAngle$};
\draw (502,191) node [anchor=north west][inner sep=0.75pt]   [align=left] {control point offsets $\mathrm{\Delta}\OptiCtrlPoint$};
\draw (367,191) node [anchor=north west][inner sep=0.75pt]   [align=left] {parameters $\OptiParameter$};
\draw  [draw opacity=0]  (37,144) -- (152,144) -- (152,165) -- (37,165) -- cycle  ;
\draw (32,148) node [anchor=north west][inner sep=0.75pt]  [font=\footnotesize] [align=left] {Optimization update};

\end{tikzpicture}
    \caption{Schematic overview for the combined parameter and shape optimization approach.}
    \label{fig:Optimization:Schematic}
\end{figure*}

A schematic overview of the optimization process is shown in \cref{fig:Optimization:Schematic}. 
The process is given for the general case of a nonlinear problem. Some simplifications can be made in the case of a linear PDE system. First, no iterative solving with the Newton-Raphson scheme is necessary. Also, as the stiffness matrices remain constant in the linear case, calculating their derivatives does not need to be repeated for different $\RotAngle$. This drastically reduces the computational effort in the linear case.

\subsubsection{Further derivatives}
Computing the analytical derivatives of $A_\mathrm{Magnet}$ with respect to $\optiVector$ from \eqref{eq:Optimization:Problem} is straightforward, because it is given as the product of the magnets' width and height, which are chosen as part of the design variables $\optiVector$. The gradient of the smoothing term $S$ does also not depend on the state equation and amounts to a sum of linear contributions.

Computing the gradients of the geometric constraints $\frac{\der\mathbf{g}}{\der\optiVector}$ is done in a general way with numerical evaluations using finite differences. Here, this is suitable, because rebuilding the geometry in IGA is fast and numerical gradients are sufficiently accurate for our simple constraints.

\section{Results}
\label{sec:Results}

\begin{figure*}[t]
    \centering
    \begin{subfigure}[t]{.44\linewidth}
      \centering
    \includegraphics[width=\linewidth]{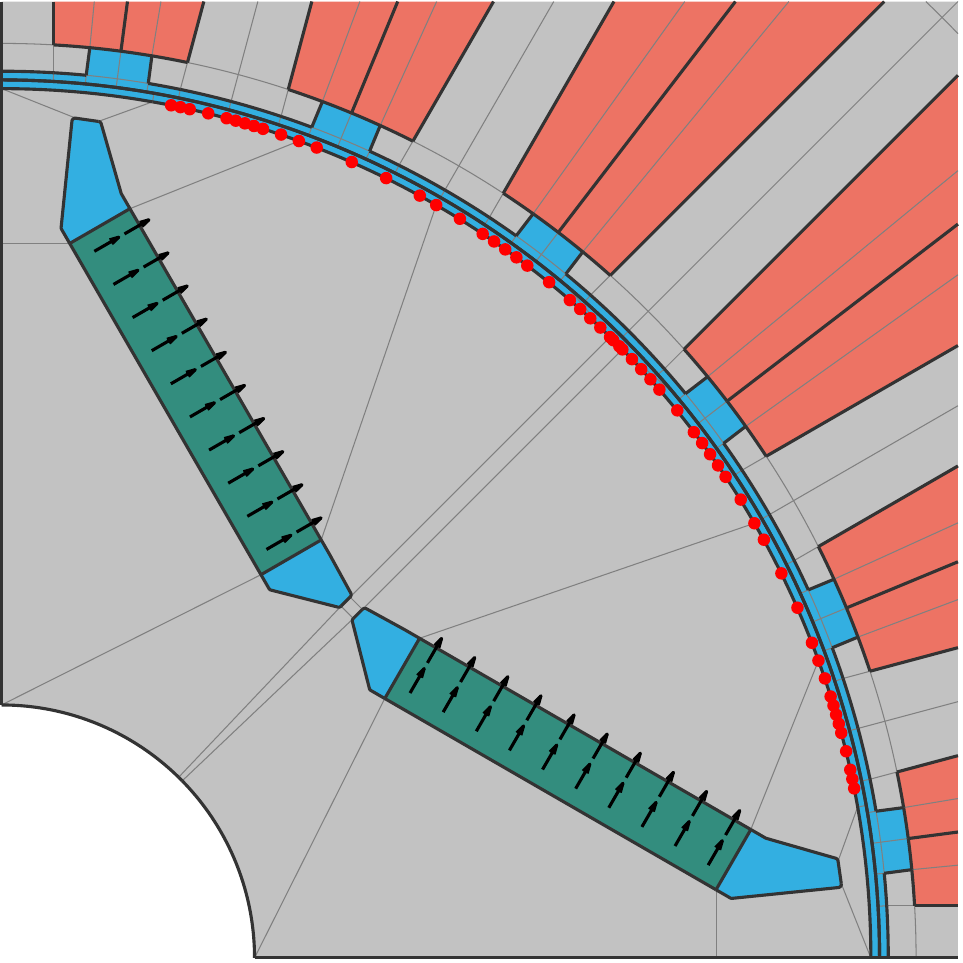}
        \caption{Patch representation of the initial motor geometry. The red points indicate the  control points for the optimization.}
      \label{fig:Results:Geometry}
    \end{subfigure} \hfill
    \begin{subfigure}[t]{.51\linewidth}
      \centering
        \includegraphics[width=\linewidth]{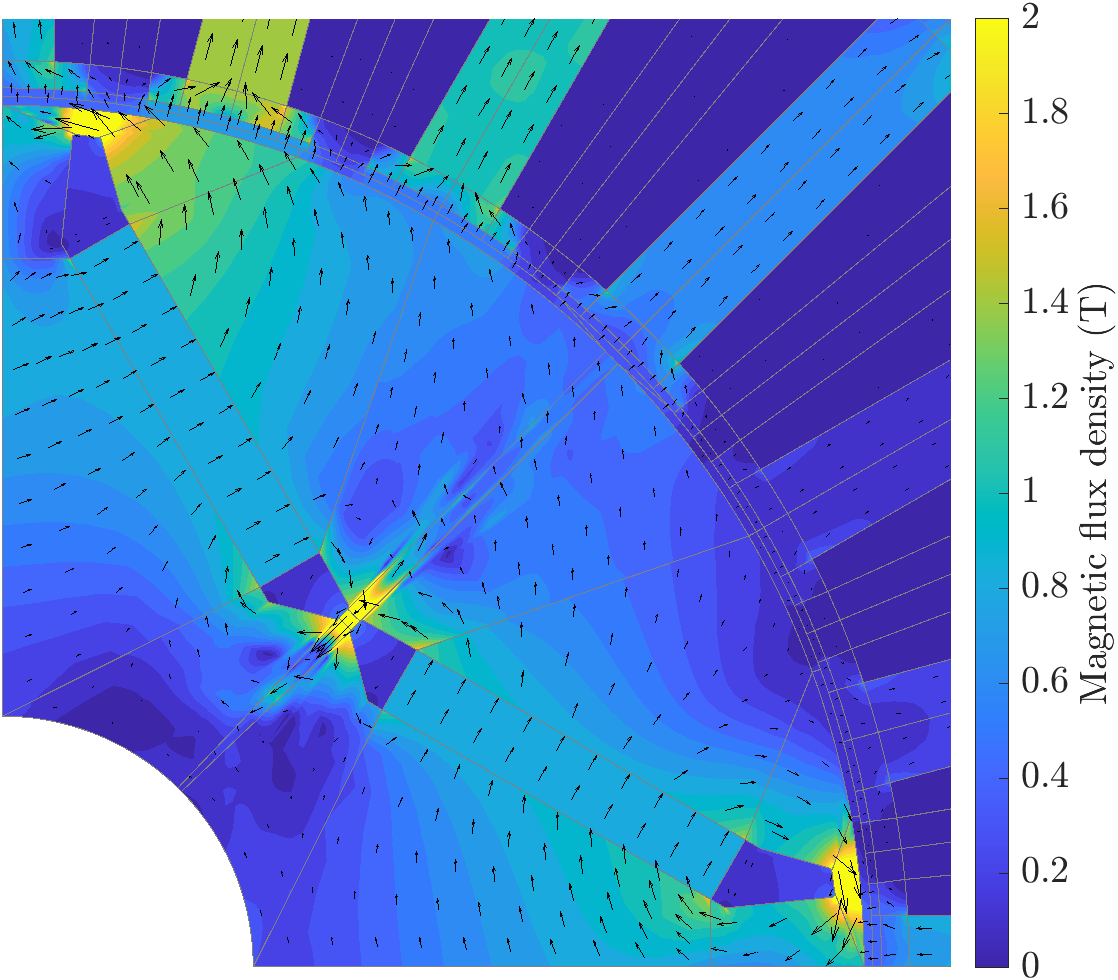}
      \caption{Magnetic flux density. The material saturates at the iron bridges, which requires the problem to be solved in a nonlinear fashion.}
      \label{fig:Results:Field}
    \end{subfigure}
    \caption{Evaluation of the initial geometry.}
    \label{fig:Results:Init}
\end{figure*}

The presented optimization workflow is implemented in \textit{MATLAB}\textsuperscript{\textregistered} using the open source package \textit{GeoPDEs} and applied on the motor geometry depicted in \cref{fig:JMAGmotor}. For more information about the implementation, the interested reader is referred to the released code provided with this publication \citep{wiesheu_2023_10160087}.

\subsection{Initial motor evaluation}

The simulations are performed with the standard M330-50A material (the equivalent M27 data from \textit{FEMM} \citep{FEMM} is used), which defines the nonlinear relative permeability $\mu_\mathrm{r}$ (inverse reluctivity) for the iron parts, and linear Nd-Fe-B magnets with $B_\mathrm{r} = \SI{1.0}{\tesla}$ and $\mu_\mathrm{r} = 1.05$. The lamination allows for disregarding the eddy currents and performing a series of magnetostatic evaluations instead. The current source density is given by an applied current of $\Iapp = \SI{3}{\ampere}$ and a winding number of $\nWind = 35$. A summary regarding the discretization and evaluation of the motor domain with IGA is given in \cref{tab:results:MotorSpecifications}.
\begin{table}[h]
    \centering
    \caption{Summary of relevant quantities regarding discretization and evaluation.}
    \label{tab:results:MotorSpecifications}
    \begin{tabular}{l|cc|}
    \cline{2-3}
                                                                & \multicolumn{1}{c|}{Rotor} & Stator \\ \hline
    \multicolumn{1}{|l|}{Number of patches}                     & \multicolumn{1}{c|}{32}    & 144    \\ \hline
    \multicolumn{1}{|l|}{Number of control points}              & \multicolumn{1}{c|}{444}   & 365    \\ \hline
    \multicolumn{1}{|l|}{Basis function degree}                 & \multicolumn{2}{c|}{2}              \\ \hline
    \multicolumn{1}{|l|}{Degrees of freedom}                    & \multicolumn{1}{c|}{4270}  & 2465   \\ \hline
    \multicolumn{1}{|l|}{Coupling degrees of freedom}           & \multicolumn{2}{c|}{52}             \\ \hline
    \end{tabular}
\end{table}

The coupling of rotor and stator domains is achieved with 26 sine and cosine functions where $n\in \mathcal{N} = \{2,6,10,...,102\}$ each. 
The rotor geometry has been further modified by introducing more control points to allow for more design freedom of the surface. The NURBS patch representation for the rotor is shown in  \cref{fig:Results:Geometry}. Additionally, the control points that are allowed to move radially during the optimization are plotted as red points. An evaluation of the magnetic flux density for a rotation angle of $\RotAngle = 0^\circ $ is shown in \cref{fig:Results:Field}.

\begin{figure}[h]
    \centering
    \input{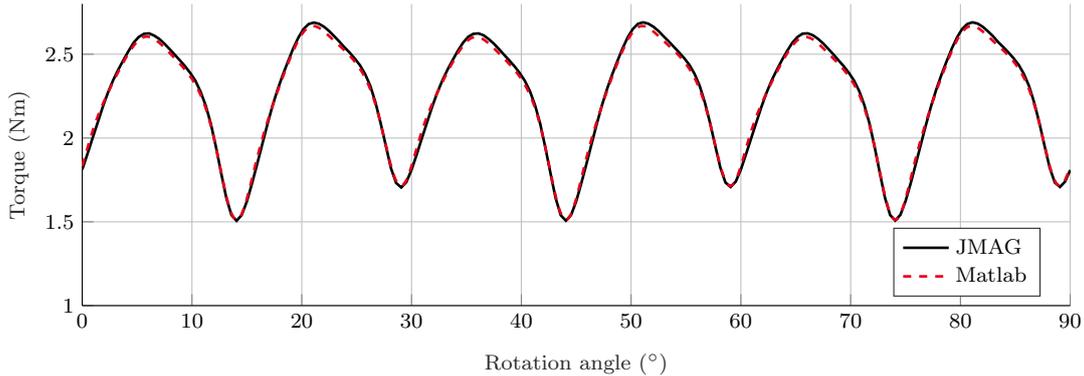}
    \caption{Comparison of the initial torque computed with Matlab and the torque computed with JMAG for the full motor geometry.}
    \label{fig:results:torque_init}
\end{figure}

The initial torque values for the full geometry calculated with \eqref{eq:Methogology:Torque} are plotted in \cref{fig:results:torque_init}. The results are validated with JMAG and show a good agreement with a mean error of $<\SI{1}{\percent}$. The machine 
has an average torque of $\MeanTorque = \SI{2.24}{Nm}$ with a quite high torque variation of $\StdTorque = \SI{0.33}{Nm}$.

Another observation is that the torque repeats every 30 degrees. This is due to the \SI{30}{\degree} symmetry of the winding pattern of the stator for synchronous operation.
Exploiting this reduces the effort for the optimization, as only a reduced number of rotation angles need to be evaluated.

\subsection{Motor optimization}
The optimization is carried out with the (semi-)analytical derivatives and gradient based optimization as explained in the previous section. Important quantities for the optimization are summarized in 
\cref{tab:results:optimization}.

\begin{table}[h]
    \centering
    \caption{Summary of the relevant quantities regarding the optimization.}
    \label{tab:results:optimization}
    \begin{tabular}{|l|c|}
    \hline
    Number of optimization parameters     & 17            \\ \hline
    Number of optimization control points & 29 (58)       \\ \hline
    Rotation angles                       & $0^\circ$, $1^\circ$, ..., $29^\circ$ \\ \hline
    $w_1$                                    & 10000       \\ \hline
    $w_2$                                    & 100        \\ \hline
    $w_3$                                    & 1000        \\ \hline
    $\Ttarget$            & \SI{2.3}{Nm}         \\ \hline
    \end{tabular}
\end{table}

\begin{figure*}[h]
    \centering
    \begin{subfigure}[t]{.24\linewidth}
      \centering
    \includegraphics[width=\linewidth]{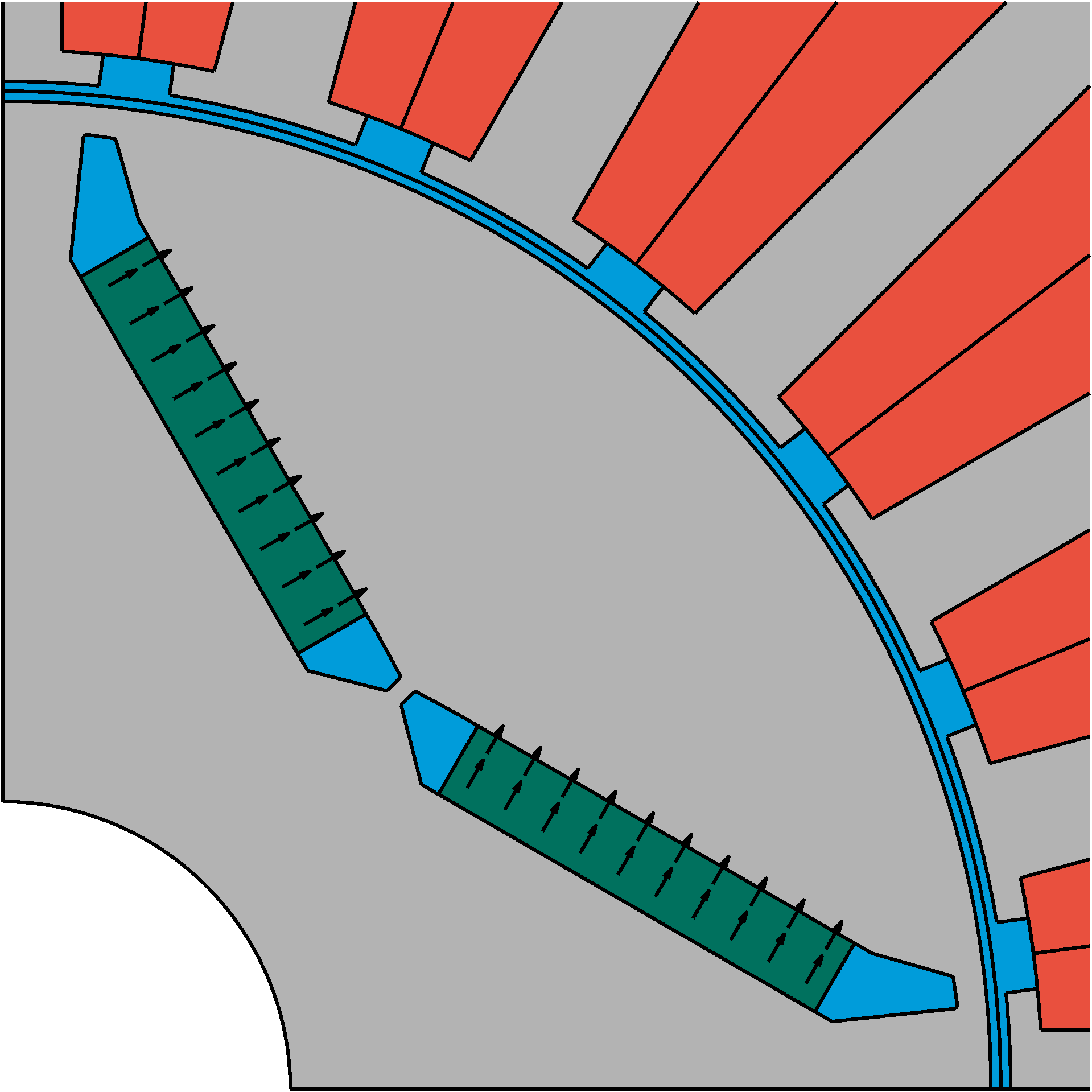}
        \caption{Iteration 0.}
      \label{fig:Results:ParamShape0}
    \end{subfigure} \hfill
    \begin{subfigure}[t]{.24\linewidth}
      \centering
        \includegraphics[width=\linewidth]{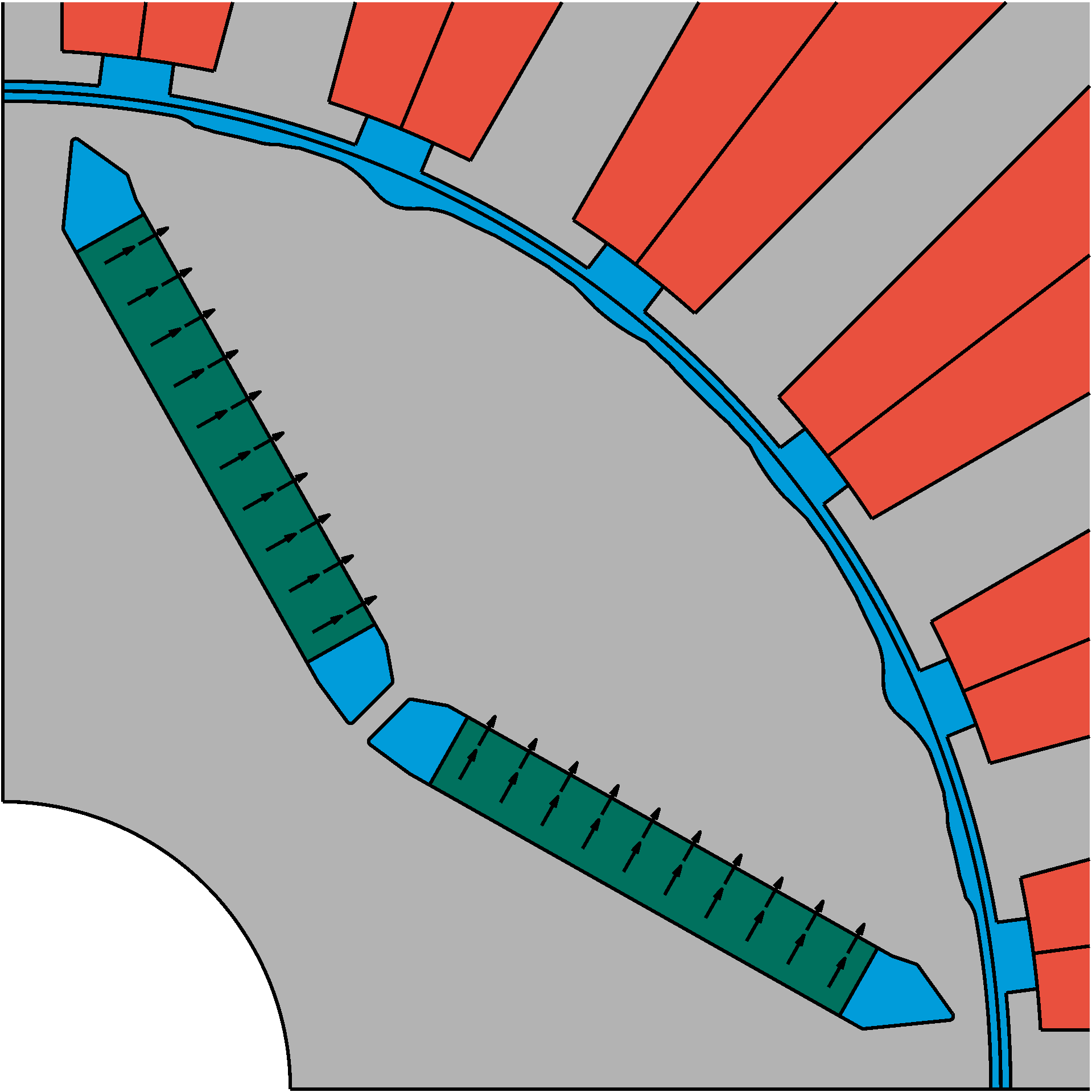}
      \caption{Iteration 5.}
      \label{fig:Results:ParamShape5}
    \end{subfigure} \hfill
    \begin{subfigure}[t]{.24\linewidth}
      \centering
      \includegraphics[width=\linewidth]{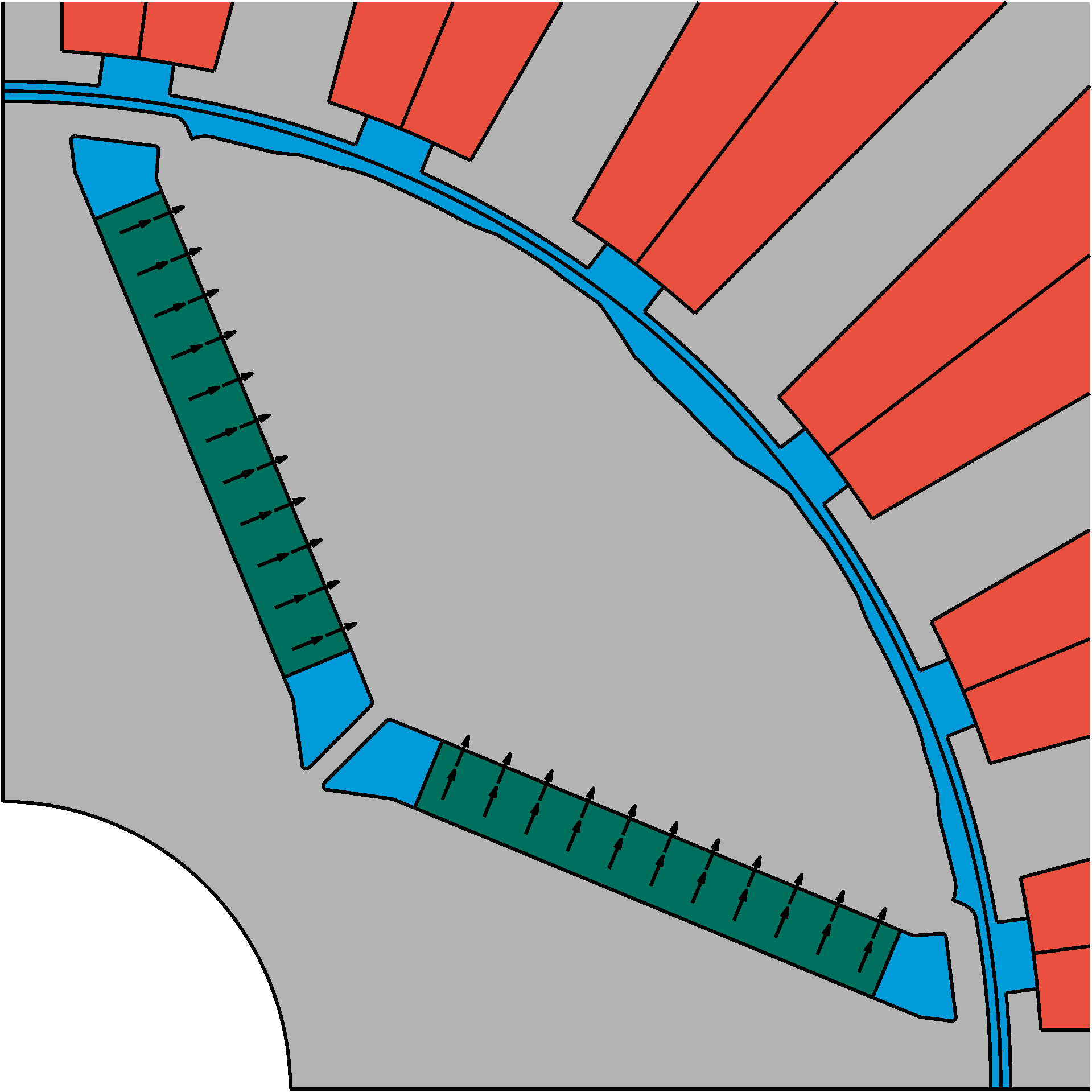}
      \caption{Iteration 50.}
      \label{fig:Results:ParamShape50}
    \end{subfigure} \hfill
    \begin{subfigure}[t]{.24\linewidth}
      \centering
      \includegraphics[width=\linewidth]{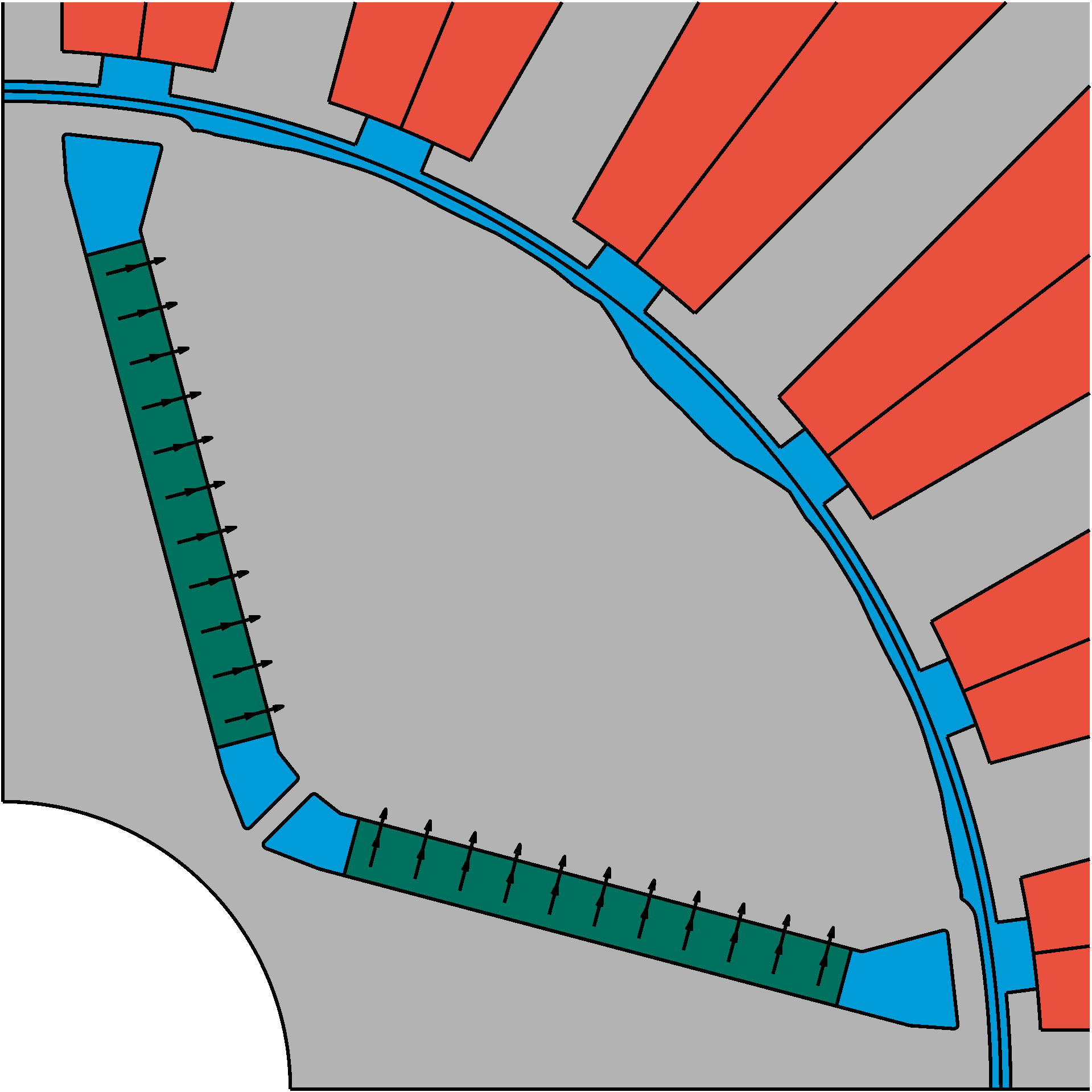}
      \caption{Iteration 200.}
      \label{fig:Results:ParamShape200}
    \end{subfigure}
    \caption{Results during the optimization process for the combined parameter and shape optimization approach.}
    \label{fig:Results:ParamShape}
\end{figure*}

There are overall 17 optimization parameters for the rotor that can be modified, e.g. those that define the permanent magnets, the air slits or the electric phase angle. A comprehensive list is given in \cref{app:tab:Parameters}. Additionally, there are 58 control points that may have a radial offset of up to \SI{1.5}{\mm}. Here we choose the geometry to be symmetric, resulting in 29 design variables for the rotor surface. Further nonlinear geometric constraints are imposed to ensure the feasibility of the design, i.e. all outer iron bridges must be at least \SI{1.5}{\mm} thick. This is achieved by nine inequality constraints restricting the air slits and magnets used for the geometry construction to remain \SI{1.5}{mm} inside of the rotor domain. The specific expressions are found in \citep{wiesheu_2023_10160087}, file \texttt{JMAG\_constraints.m}. The weights $w_1$, $w_2$ and $w_3$ are chosen such that the objective function and constraints are in a comparable order of magnitude. This facilitates the numerical stability of the gradient based solver. The torque is evaluated at each degree, so that the rotation angles are set to $\beta \in \{0, 1, ..., 29\}$ due to the torque's periodicity. The mean torque shall be improved from the initial $\MeanTorque= \SI{2.24}{Nm}$ 
to $\Ttarget = \SI{2.3}{Nm}$.

\begin{figure*}[h]
    \centering
    \begin{subfigure}[t]{.24\linewidth}
      \centering
    \includegraphics[width=\linewidth]{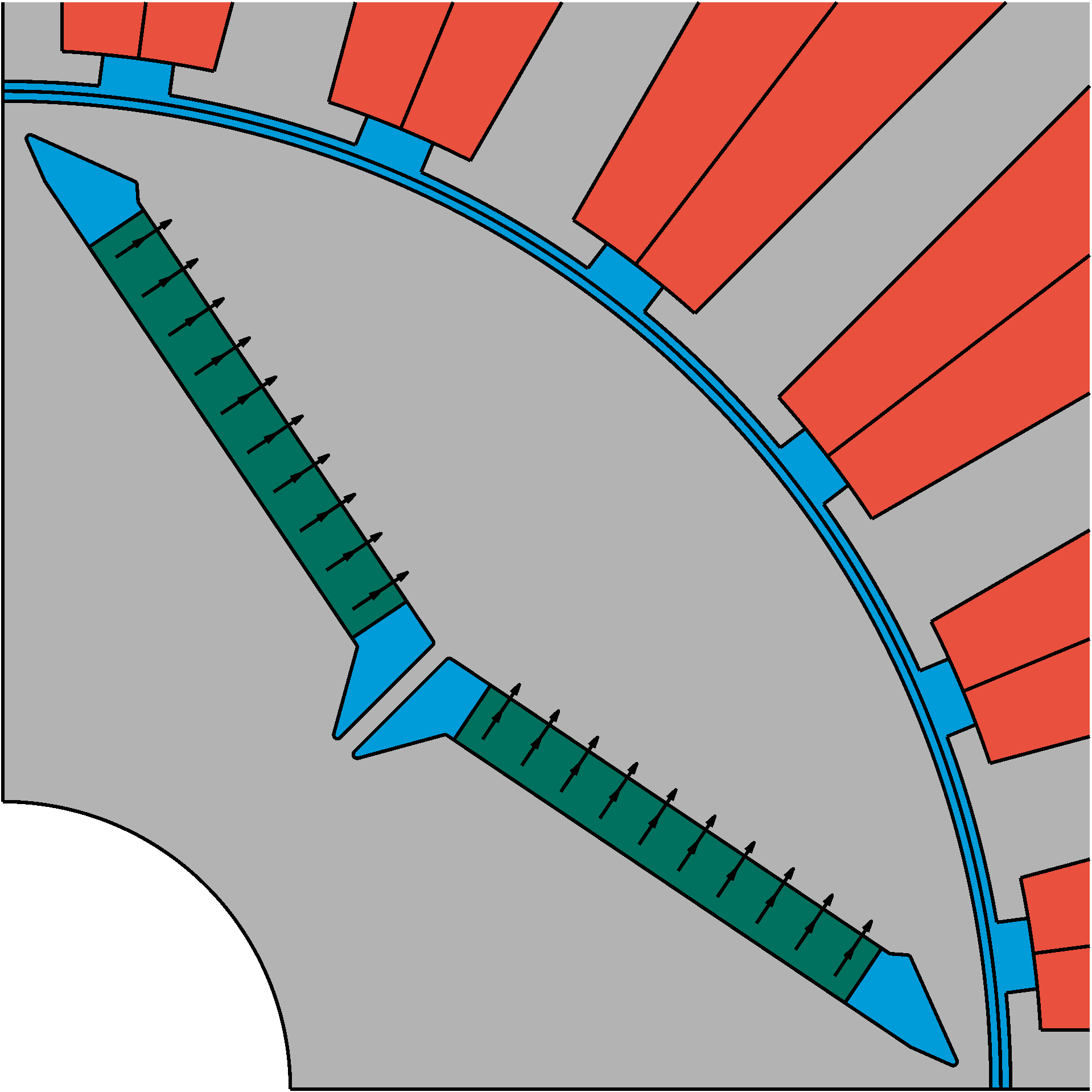}
        \caption{Parameter Optimized}
      \label{fig:Results:ParamOpt}
    \end{subfigure} \hfill
    \begin{subfigure}[t]{.24\linewidth}
      \centering
        \includegraphics[width=\linewidth]{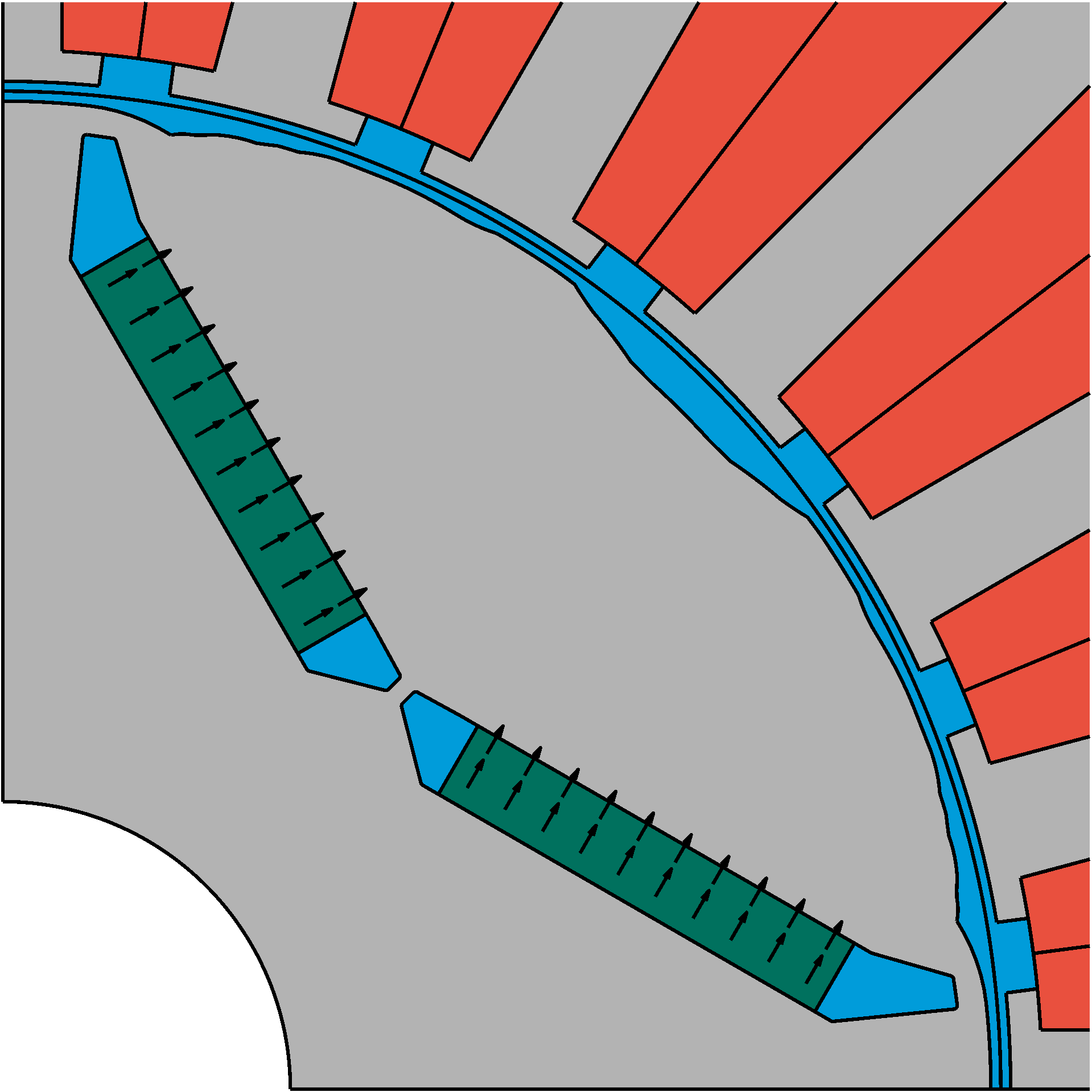}
      \caption{Shape Optimized}
      \label{fig:Results:ShapeOpt}
    \end{subfigure} \hfill
    \begin{subfigure}[t]{.24\linewidth}
      \centering
      \includegraphics[width=\linewidth]{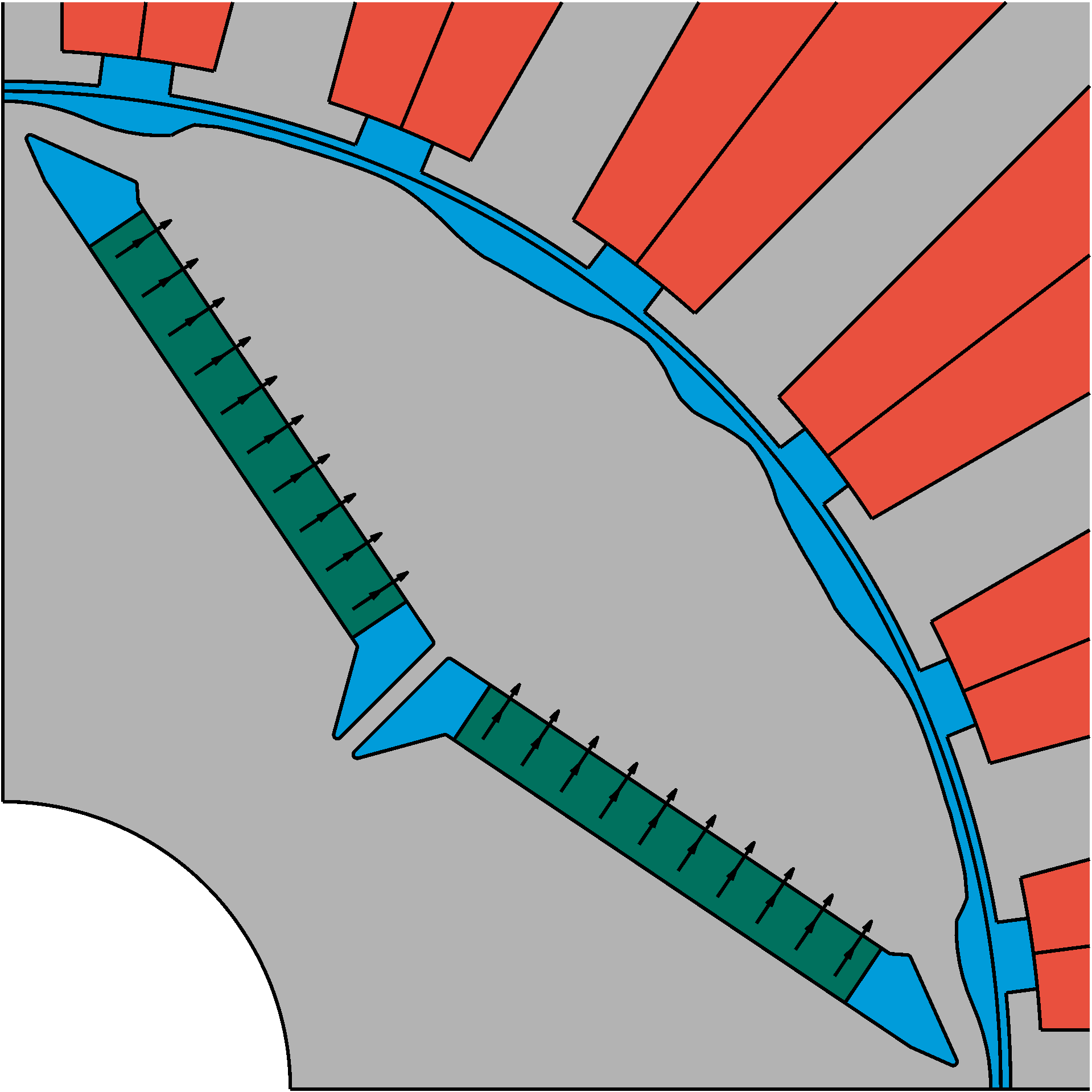}
      \caption{Parameter then shape optimized}
      \label{fig:Results:ParamThenShapeOpt}
    \end{subfigure} \hfill
    \begin{subfigure}[t]{.24\linewidth}
      \centering
      \includegraphics[width=\linewidth]{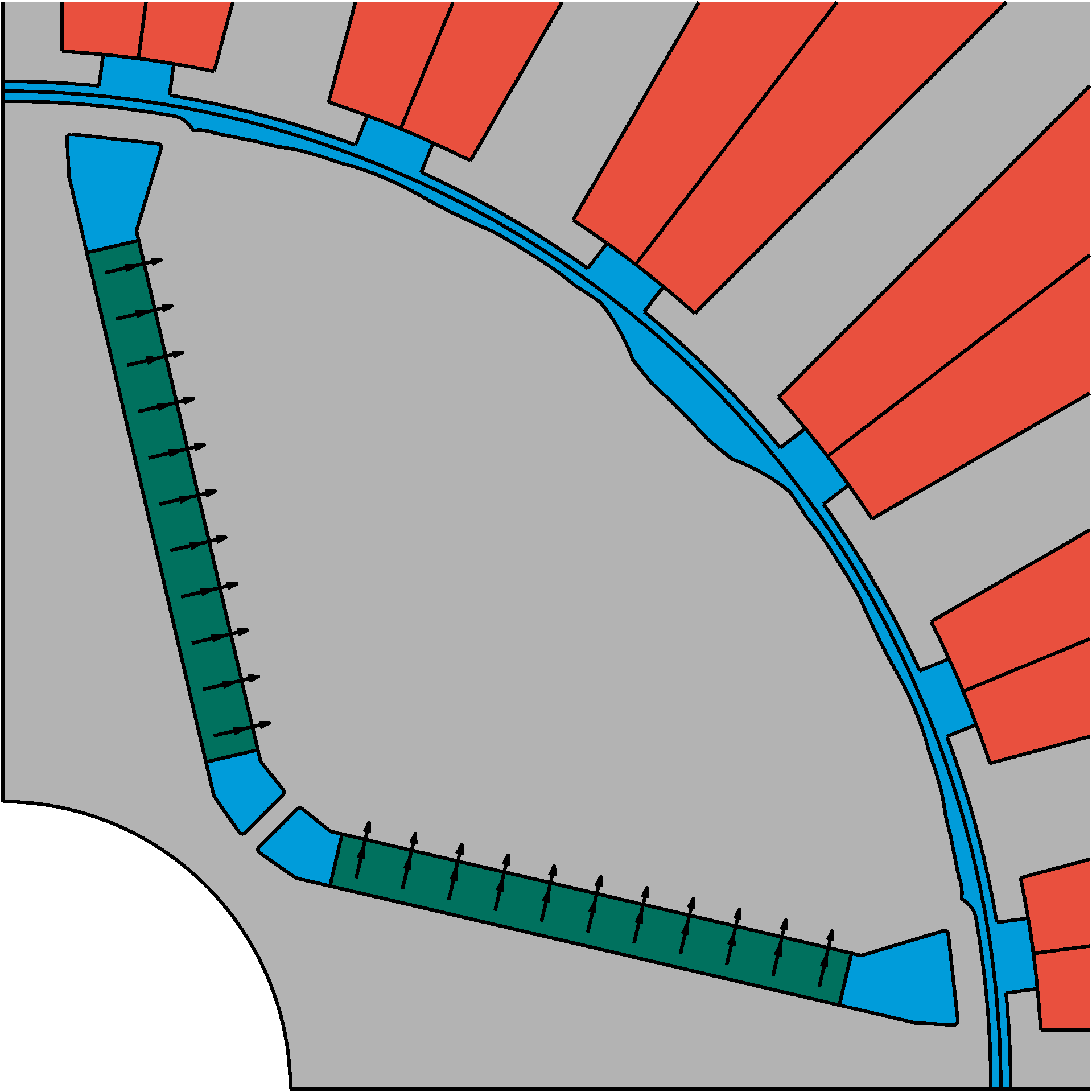}
      \caption{Parameter and shape optimized}
      \label{fig:Results:ParamShapeOpt}
    \end{subfigure}
    \caption{Comparison of final geometries for the different optimization methods.}
    \label{fig:Results:OptComparison}
\end{figure*}

The optimization process is carried out on a 8-core laptop (Intel\textsuperscript{\textregistered} Core\textsuperscript{\texttrademark} i7-1165G7@2.8GHz) with \SI{16}{GB} RAM in MATLAB R2024a. As optimization method, the \textit{interior-point} solver from \textit{MATLAB}\textsuperscript{\textregistered} is used. The optimization function, constraints and their gradients are provided to the \textit{fmincon} solver. Some intermediate geometries during the optimization process for the combined parameter and shape approach are presented in \cref{fig:Results:ParamShape}. The optimization starts with the initial motor geometry given by the standard V-shape template (\cref{fig:Results:ParamShape0}). It can be seen that multiple geometry features, such as the magnet position, air pockets and rotor surface, are changed simultaneously. After 200 iterations the optimization converges to a design with larger air slits, thinner magnets and adjusted rotor surface (\cref{fig:Results:ParamShape200}).

\subsection{Comparison}
For further comparison with conventional optimization methods, the optimization process is performed with several other configurations, i.e. solely parameter optimization, solely shape optimization and sequential parameter and shape optimization. The final results for the different optimization methods are summarized in \cref{tab:results:comparison}. Note that the optimizer could also converge to a different optimum for other initial conditions, because of the nature of gradient-based optimization methods.

\begin{table*}[h]
\caption{Optimization results for the studied approaches. The values are given for the reduced motor geometry (only one quarter).}
\label{tab:results:comparison}
\centering
\begin{tabular}{|l||c|c|c||c|c|c|c|}
\hline
Method & $\fopt$ & $A_\mathrm{Magnet} (\unit{m^2})$ & $\StdTorque (\unit{Nm})$  & $\MeanTorque (\unit{Nm})$  & Iterations & Evaluations & Time \\
\hline
\hline
Initial & 10.124 & 1.760e-04 & 0.0836 & 0.561 & - & - & - \\
\hline
Parameter & 4.081 & 1.591e-04 & 0.0249 & 0.575 & 51 & 65 & \SI{35}{min} \\
\hline
Shape & 2.331 & 1.760e-04 & 0.0057 & 0.501 & 111 & 190 & \SI{2}{h} \SI{7}{min} \\
\hline
Sequential parameter and shape & 1.756 & 1.591e-04 & 0.0008 & 0.531 & 51 + 220 & 65 + 353 & \SI{35}{min} + \SI{3}{h} \SI{40}{min} \\
\hline
Combined parameter and shape & 1.596 & 1.428e-04 & 0.0013 & 0.576 & 284 & 429 & \SI{4}{h} \SI{16}{min} \\
\hline
JMAG EA parameter & 4.509 & 1.747e-04 & 0.0276 & 0.576 & 2500 & 2500 & \SI{14}{h}\footnotemark \\
\hline
\end{tabular}
\end{table*}
\footnotetext[1]{Run on a slower working station due to high memory usage.}
The final geometries for the different optimization methods are shown in \cref{fig:Results:OptComparison} and explained in the following:
\begin{itemize}
    \item The parameter optimized design achieves a reduction of $\Amagnet$ by \SI{9.6}{\percent} and a reduction of $\StdTorque$ by \SI{70.2}{\percent} within 51 iterations. This optimization is fastest with a calculation time of about 35 minutes. The result is shown in \cref{fig:Results:ParamOpt}.
    \item The shape optimized design (\cref{fig:Results:ShapeOpt}) reaches an even higher $\StdTorque$-reduction of \SI{94.2}{\percent} by modifications of the rotor surface. Since the magnet mass remains unchanged, the requirements on the average torque are relaxed to $\Ttarget = \SI{2}{Nm}$. This is necessary, as modifications in the rotor shape alone do not allow for significant increases of $\MeanTorque$, so the design would be infeasible. One parameter that may be optimized here is the electric phase angle $\CurrentAngle$ as this value is independent of the geometry. 
    \item The sequential parameter-then-shape optimized design (\cref{fig:Results:ParamThenShapeOpt}) achieves a reduction of $\Amagnet$ by \SI{9.6}{\percent} and a even greater reduction of $\StdTorque$ by \SI{98.4}{\percent}. But again, we set the reduced $\Ttarget = \SI{2}{Nm}$, as the shape optimization starts from a local minimum given by the previous optimization. Therefore it is necessary to decrease $\Ttarget$ in order to allow for further smoothing of the torque profile. Also $\CurrentAngle$ may be changed again during the shape optimization.
    \item The best optimization result is given by the combined parameter-and-shape optimized design (\cref{fig:Results:ParamShapeOpt}) with the methodology presented in the previous chapter. Due to the increased solution space, this optimization process needs the most function evaluations and also takes longest with a simulation time of more than 4 hours. Overall, about \SI{7}{\percent} of the time is needed for reevaluations of the geometry and basis functions, \SI{65}{\percent} for rebuilding of matrices and repeated solution of the state equation at every angle, \SI{3}{\percent} for the solution of the adjoint system, \SI{23}{\percent} for the calculation of the objective function and required derivatives and \SI{2}{\percent} for miscellaneous tasks (such as updating plots and storing data). The simulations indicate that even higher reductions of $\Amagnet$ by \SI{18.9}{\percent} and $\StdTorque$ by \SI{99.0}{\percent} are possible. This is again a major improvement compared to the previous methods.
    \item For comparison, a simple evolutionary algorithm (EA) optimization for the parameters has been performed with JMAG. The commercial tool reaches reductions of $\Amagnet$ by \SI{0.7}{\percent} and $\StdTorque$ by \SI{67}{\percent} within 14 hours. Note that these values are given for a rough comparison and should be treated with care, as setting the objective function and constraints is different in JMAG. The EA also finds multiple designs and results could be further improved with more computation time.
\end{itemize}

Note that many processes in the Matlab prototype code are not parallelized (such as assembling matrices and derivatives or solving the problem for different $\RotAngle$). It is expected that the computation time can be drastically reduced to mere hours or even minutes if the code is rewritten with performance in mind. This is also reflected in the distribution of the simulation times. Most time is needed for the rebuilding of matrices, whereas the solving of linear systems -- which is optimized in Matlab -- is almost negligible. It is also possible to speed up the optimization process by reducing the amount of angles that are evaluated, e.g. by only evaluating every third degree. However, this can introduce higher order harmonics in the torque which are not covered by the sampled angles, resulting in an underestimation of $\StdTorque$. 

\begin{figure}[h]
    \centering

%
%
\definecolor{mycolor1}{rgb}{0.36471,0.52157,0.76471}%
\definecolor{mycolor2}{rgb}{0.31373,0.71373,0.58431}%
\definecolor{mycolor3}{rgb}{0.97255,0.72941,0.23529}%
\definecolor{mycolor4}{rgb}{0.91373,0.31373,0.24314}%
\definecolor{mycolor5}{rgb}{0.30100,0.74500,0.93300}%
\definecolor{mycolor6}{rgb}{0.63500,0.07800,0.18400}%
\definecolor{mycolor7}{rgb}{0.00000,0.44700,0.74100}%
\definecolor{mycolor8}{rgb}{0.85000,0.32500,0.09800}%
\begin{tikzpicture}

\begin{axis}[%
width=14cm,
height=6cm,
at={(1.137in,0.814in)},
scale only axis,
xmin=0,
xmax=30,
xlabel style={font=\color{white!15!black}},
xlabel={$\text{Rotation angle (}^\circ\text{)}$},
ymin=0.5,
ymax=2.7,
ylabel style={yshift=-0.4cm, font=\color{white!15!black}},
ylabel={Torque (Nm)},
axis background/.style={fill=white},
axis x line*=bottom,
axis y line*=left,
xmajorgrids,
ymajorgrids,
legend style={at={(0.97,0.03)}, anchor=south east, legend cell align=left, align=left, draw=white!15!black}
]
\addplot [color=black, line width=1.5pt]
table[x index=0,y index=1,col sep=comma] {plots/allData.csv};
\addlegendentry{Initial}

\addplot [color=mycolor1, line width=1.5pt,legend image post style={mark={*}}]
table[x index=0,y index=2,col sep=comma] {plots/allData.csv};

\addplot[only marks, mark=*, mark options={}, mark size=2pt, color=mycolor1, fill=mycolor1, forget plot] table[x index=0,y index=1,col sep=comma] {plots/optData.csv};
\addlegendentry{Parameter optimized}

\addplot [color=mycolor2, line width=1.5pt,legend image post style={mark={square*}}]
  table[x index=0,y index=3,col sep=comma] {plots/allData.csv};
\addplot[only marks, mark=square*, mark options={}, mark size=2pt, color=mycolor2, fill=mycolor2, forget plot] table[x index=0,y index=2,col sep=comma] {plots/optData.csv};
\addlegendentry{Shape optimized}

\addplot [color=mycolor3, line width=1.5pt,legend image post style={mark={triangle*}}]
  table[x index=0,y index=4,col sep=comma] {plots/allData.csv};
\addplot[only marks, mark=triangle*, mark options={}, mark size=2pt, color=mycolor3, fill=mycolor3, forget plot] table[x index=0,y index=3,col sep=comma] {plots/optData.csv};
\addlegendentry{Parameter then shape optimized}

\addplot [color=mycolor4, line width=1.5pt, legend image post style={mark={diamond*}}]
  table[x index=0,y index=5,col sep=comma] {plots/allData.csv};
\addplot[only marks, mark=diamond*, mark options={}, mark size=2pt, color=mycolor4, fill=mycolor4, forget plot] table[x index=0,y index=4,col sep=comma] {plots/optData.csv};
\addlegendentry{Parameter and shape optimized}

\end{axis}

\end{tikzpicture}
    \caption{Comparison of the torque profiles for the final geometries obtained by the different optimization methods.}
    \label{fig:results:torque_comparison}
\end{figure}
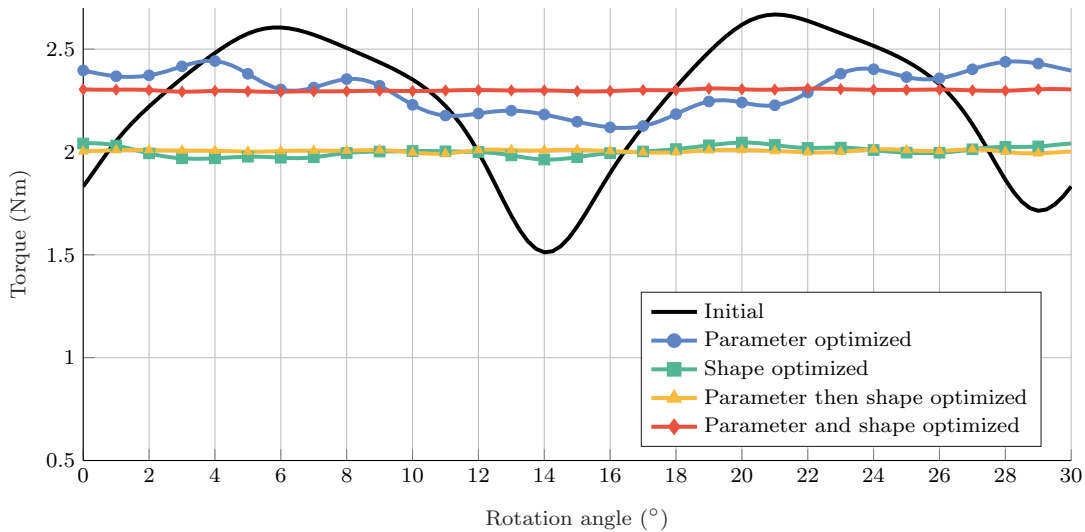

Finally, the torque curves for the resulting geometries are given in \cref{fig:results:torque_comparison}. The markers indicate the rotation angles, which are used for evaluation during the optimization. 
One can observe the different constraints of $\Ttarget\geq\SI{2.3}{Nm}$ and $\Ttarget\geq\SI{2.0}{Nm}$ for the different optimization methods. As explained, reducing the target torque was necessary for the shape optimization, because only the surface may be changed in this optimization step. If the optimization includes changes of the rotor shape, a significant reduction of the torque ripple is achievable compared to the parameter optimized geometry. If parameter and shape optimization is combined in some way, the torque ripple can be reduced to almost zero in simulations.

Note that our optimization results may not be optimal from different physical points of view, such as mechanical and thermal stresses or manufacturability. The goal of this paper is to provide insights into a new optimization procedure and its application to motor design. 
Including further relevant quantities, such as mechanical stresses, with a multiphysical optimization is left for future work.

\subsection{Prospect for further work}
The presented results show a novel spline-based combined parameter/shape approach applied to the magnetic field in a PMSM. The results are already very promising and there are various prospects for further research.

One straightforward idea is to also include the stator geometry in all optimization steps.
This would also allow to introduce a scaling factor for the radial dimension and thus enable different motor sizes and a reduction of the total motor cost. By further defining a relationship between axial machine length and torque (currently a linear one is used) to account for fringing effects, all relevant geometric parameters can be determined in one optimization process. So far, symmetry of the rotor has been assumed, but the motor can also be designed in non-symmetrical ways to achieve new, out-of-the-box designs. Further, robust optimization can be studied, which contributes to easier manufacturability by allowing larger tolerances. Another idea is, to include second order derivatives to increase the convergence rate and enable efficient multi-objective optimization.

One necessary step in achieving smaller motor geometries, higher power densities, and pushing designs to their physical limits is to include the calculation of mechanical stresses and impose restrictions on the maximum permissible values. 
Since the optimization method is not limited to magnetism, incorporating mechanical stress constraints is possible in the same way as described. Calculating stress derivatives with respect to NURBS control points is already known from the literature \citep{Pal2021xx}. Overall, the presented method offers numerous advantages, and there are several ways to further speed up the computations, extend the method, or apply it to different fields in engineering.
\section{Conclusions}
\label{sec:Conclusions}
This paper presented a novel and efficient optimization procedure, which takes advantage of the intrinsic properties of Isogeometric Analysis. IGA allows for straightforward geometry modifications by adjusting control points, eliminates the need for meshing and enables an exact representation of curved segments, which is particularly useful for rotating machines. 

The method described shows a way to efficiently combine parameter and shape optimization with gradient based optimization. The optimization process is exemplary applied to a PMSM with the goal to minimize the usage of permanent magnets and reduce the torque ripple while maintaining a specified target torque. It has been shown that the combination of parameter and shape optimization yields the best possible results. Also an extension of the method and applications in different fields of physics are possible.

\section*{Declarations}

\subsection*{Ethics approval and consent to participate:}
Not applicable
\subsection*{Consent for publication:}
Not applicable
\subsection*{Availability of data and materials:}
All data and materials are publicly available via Github/Zenodo \citep{wiesheu_2023_10160087}. 

\subsection*{Authors' contributions:}
\begin{tabular}{p{0.22\linewidth} p{0.75\linewidth}}
         \textbf{Michael Wiesheu}: & Methodology, Software, Formal analysis, Investigation, Data Curation, Writing - Original Draft, Visualization \\
        \textbf{Theodor Komann}: & Methodology, Software, Validation, Formal analysis, Investigation, Data Curation, Writing - Original Draft, Visualization \\
        \textbf{Melina Merkel}: & Methodology, Software, Writing - Review \& Editing, Funding acquisition \\
        \textbf{Sebastian Schöps}: & Conceptualization, Resources, Writing - Review \& Editing, Supervision, Funding acquisition \\
        \textbf{Stefan Ulbrich}: & Conceptualization, Formal analysis, Resources, Writing - Review \& Editing, Supervision, Funding acquisition \\
        \textbf{Idoia Cortes Garcia}: & Conceptualization, Writing - Review \& Editing, Supervision 
\end{tabular}

\subsection*{Competing interests:}
The authors declare that they have no competing interests.

\subsection*{Funding}
The work of Michael Wiesheu and Theodor Komann is partially supported by the joint DFG/FWF Collaborative Research Centre CREATOR (DFG: Project-ID 492661287/TRR 361; FWF: 10.55776/F90) at TU Darmstadt, TU Graz and JKU Linz.
Further support is given by the Graduate School CE within the Centre for Computational Engineering at TU Darmstadt.

\subsection*{Acknowledgements}
The authors thank JMAG for the permission to use the machine model and the cooperational software licences.

\appendix
\section{Derivatives}
\label{app:appendix}


\subsection{Jacobian for Newton-Raphson}
\label{app:NewtonJacobian}
The derivative of the stiffness matrix $\stiffnessSymbol_{ij}$ is given by
\begin{equation}
    \frac{\operatorname{d} K_{ij}}{\operatorname{d} u_{k}} =\frac{\operatorname{d}}{\operatorname{d} u_k}\int _{\Omega } \nu ( B) \nabla N_{i} \cdotp \nabla N_{j}\operatorname{d} \Omega =\int _{\Omega }\underbrace{\frac{\partial \nu ( B)}{\partial B}}_{BH-curve}\frac{\operatorname{d} B}{\operatorname{d} u_{k}} \nabla N_{i} \cdotp \nabla N_{j}\operatorname{d} \Omega 
\end{equation}
where $B=B( \solutionVector)$ and
\begin{align}
    \frac{\operatorname{d} B}{\operatorname{d} u_{k}} & =\frac{\left(\sum _{l} \nabla N_{l} u_{l}\right) \cdotp \nabla N_{k}}{B}, & B= & \sqrt{\left(\sum N_{l,x} u_{l}\right)^{2} +\left(\sum N_{l,y} u_{l}\right)^{2}}.
\end{align}

\subsection{Sensitivities of stiffness matrix and right-hand side}
\label{app:sensitivityStiffnessMat}

The derivatives of the stiffness matrix are given by the product rule, which results in four different contributions:
\begin{equation}
    \frac{\operatorname{d} K_{ij}}{\operatorname{d} C_{kd}} =K_{ij,kd}^{( 1)} +K_{ij,kd}^{( 2)} +K_{ij,kd}^{( 3)} +K_{ij,kd}^{( 4)} \label{eq:app:ProductRule}
\end{equation}

The four components of \eqref{eq:app:ProductRule} are calculated as: 
\begin{align}
K_{ij,kd}^{( 1)} & =\int _{\hat{\Omega }} -\nu (\mathbf{D}_{kd} \nabla N_{i}) \cdotp \nabla N_{j} \ |\mathbf{J}_{F} |\ \operatorname{d}\hat{\Omega }\\
K_{ij,kd}^{( 2)} & =\int _{\hat{\Omega }} -\nu \nabla N_{i} \cdotp (\mathbf{D}_{kd} \nabla N_{j}) \ |\mathbf{J}_{F} |\ \operatorname{d}\hat{\Omega }\\
K_{ij,kd}^{( 3)} & =\int _{\hat{\Omega }} \nu \nabla N_{i} \cdotp \nabla N_{j} \ |\mathbf{J}_{F} |\ \mathrm{Tr}(\mathbf{D}_{kd}) \ \operatorname{d}\hat{\Omega }\\
K_{ij,kd}^{( 4)} & =\int _{\hat{\Omega }}\underbrace{\frac{\partial \nu }{\partial B}}_{BH-curve}\frac{dB}{dC_{kd}} \nabla N_{i} \cdotp \nabla N_{j} \ |\mathbf{J}_{F} |\ \operatorname{d}\hat{\Omega }
\end{align}
with 
\begin{equation}
    \mathbf{D}_{kd} =\begin{cases}
    \begin{pmatrix}
    \nabla G_{k} & \mathbf{0}
    \end{pmatrix} & d=x\\
    \begin{pmatrix}
    \mathbf{0} & \nabla G_{k}
    \end{pmatrix} & d=y
    \end{cases}
\end{equation}
and
\begin{align}
\frac{\operatorname{d} B}{\operatorname{d} C_{kd}} & =\frac{\operatorname{d}}{\operatorname{d} C_{kd}}\sqrt{\left(\sum\nolimits _{l} N_{l,x} u_{l}\right)^{2} +\left(\sum\nolimits _{l} N_{l,y} u_{l}\right)^{2}}\\
 & =\frac{\left(\sum _{l} \nabla N_{l} u_{l}\right) \cdotp \left(\sum _{l} -\mathbf{D}_{kd} \nabla N_{l} u_{l}\right)}{\sqrt{\left(\sum _{l} N_{l,x} u_{l}\right)^{2} +\left(\sum _{l} N_{l,y} u_{l}\right)^{2}}}
\end{align}

In analogous fashion, the derivative of the right hand side with respect to the control points are given by
\begin{equation}
    \frac{\operatorname{d} b_{i}}{\operatorname{d} C_{kd}} =b_{i,kd}^{( 1)} +b_{i,kd}^{( 2)}
\end{equation}
where $b_{i,kd}^{(1)}$ and $b_{i,kd}^{(2)}$ can be calculated similarly.  It follows 
\begin{align}
b_{i,kd}^{( 1)} = & \nu \int _{\hat{\Omega }} -(\mathbf{D}_{kd} \nabla N_{i}) \cdotp \mathbf{B}_{\mathrm{r}}^{\bot } |\mathbf{J}_{F} |\ \operatorname{d}\hat{\Omega }\\
b_{i,kd}^{( 2)} = & \nu \int _{\hat{\Omega }} \nabla N_{i} \cdotp \mathbf{B}_{\mathrm{r}}^{\bot } |\mathbf{J}_{F} |\ \mathrm{Tr}(\mathbf{D}_{kd}) \ \operatorname{d}\hat{\Omega }
\end{align}
In order to address changes of the magnetization direction, the derivatives of \eqref{eq:methodology:brt} with respect
to the magnetization angle must be considered. This is given by
\begin{equation}
    \frac{\operatorname{d} b_{\rt,i}}{\operatorname{d} \alpha } =B_{\mathrm{r}} \nu \int _{\Omega } \nabla N_{i} \cdotp \begin{pmatrix}
    -\cos( \alpha )\\
    -\sin( \alpha )
    \end{pmatrix}\operatorname{d} \Omega .
\end{equation}
Finally, changes of the electric phase angle $\CurrentAngle$ are considered by the derivatives of \eqref{eq:Methodology:bst}:
\begin{equation}
    \frac{\operatorname{d} b_{\st,i}}{\operatorname{d} \CurrentAngle }=\int _{\Omega } \frac{\Iapp \nWind}{\Acoil}\cos\left(\polePair\RotAngle +\CurrentAngle +\frac{2\pi }{3} k\right) N_{i}\operatorname{d} \Omega
\end{equation}

\subsection{Motor definitions}
\begin{table}[h]
    \centering
    \caption{Material definitions.}
    \begin{tabular}{|l|c|c|c|}
    \hline
    Material & $\mu_\mathrm{r}$                     & $B_\mathrm{r}$ (T) & Color \\ \hline
    Iron     & nonlinear (M27) & - & Gray  \\ \hline
    Magnet   & 1.05                      & 1.0  & Green \\ \hline
    Copper   & 1.0                       & -    & Red   \\ \hline
    Air      & 1.0                       & -    & Blue  \\ \hline
    \end{tabular}
     \label{app:tab:Materials}
\end{table}

\begin{table}[h]
    \centering
    \caption{Initial motor parameters and bounds for the optimization. Lengths are given in \SI{}{mm}, angles in degrees.}
    \begin{tabular}{|l|c|c|c|}
    \hline
    Parameter name & initial & min & max \\
    \hline
    DMAG & 30 & 20 & 40 \\
    DSLIT5 & 1 & 0.5 &  6 \\
    DSLIT6 & 2 & 1.5 &  8 \\
    LSLIT1 & 6.4 & 1 & 8 \\
    LSLIT2 & 4.3 & 1 & 8\\
    MA & 150 & 100 & 160 \\
    MT1 & 4 & 1.5 & 10 \\
    MW1 & 22 & 10 & 27 \\
    OPERATING\_ANGLE  & 0 & -20 & 20 \\
    RA1 & 144 & 130 & 170 \\
    RA2 & 166 & 150 & 179 \\
    RS  & 1 & 1 & 4 \\
    RW2 & 1 & 0.5 & 4 \\
    RW3 & 1 & 0.5 & 4 \\
    RW4 & 1 & 0.5 & 2 \\
    RW5 & 1 & 0.5 & 2 \\
    WMAG & 4 & 3 & 8 \\
    \hline
    LENGTH & 35 &  &  \\
    RD1 & 100 & & \\
    RD2 & 29.1 & & \\
    RF & 0.2 &  &  \\
    SD1 & 204 & & \\
    SD2 & 102 & & \\
    ST & 1.6 & & \\
    SW1 & 6 & & \\
    SW2 & 3.6 & & \\
    SW4 & 25.5 & & \\
    \hline
    \end{tabular}
    \label{app:tab:Parameters}
\end{table}

\end{document}